\documentclass[11pt, a4paper, oneside, reqno]{amsart}
\usepackage{amsmath,amsfonts, graphicx,fancyhdr,mdframed}
\usepackage{algorithmic}
\usepackage{algorithm}
\usepackage{caption}
\usepackage{array}
\usepackage[caption=false,font=normalsize,labelfont=sf,textfont=sf]{subfig}
\usepackage{textcomp}
\usepackage{stfloats}
\usepackage{url}
\usepackage{verbatim}
\usepackage{graphicx}
\usepackage[backend=biber,style=ieee,sorting=none]{biblatex}
\title[]{Functional Architecture and Implementation of an Autonomous Emergency Steering System}

\author{Chris van der Ploeg}
\author{Ron Wouters}
\author{Mart Baars} 
\author{Anh-Lam Do}
\author{François Fauvel}
\thanks{C.J. van der Ploeg is with the Department of Mechanical Engineering, Eindhoven University of Technology,
({\tt C.J.v.d.Ploeg$@$tue.nl}) and the Netherlands Organisation for Applied Scientific Research, Integrated Vehicle Safety Group, 5700 AT Helmond, The Netherlands. ({\tt chris.vanderploeg@tno.nl}),
M.C.W. Baars and R.M.T. Wouters are with the Netherlands Organisation for Applied Scientific Research, Integrated Vehicle Safety Group, 5700 AT Helmond, The Netherlands. ({\tt mart.baars@tno.nl, ron.wouters@tno.nl}),
A.L. Do and F. Fauvel are with Renault S.A.S, 1 av. du Golf, 78288 Guyancourt, France. ({\tt anh-lam.do@renault.com, francois.fauvel@renault.com})}
\begin{document}
\maketitle

\begin{abstract}
    Autonomous emergency steering (AES) systems have the promising potential to further reduce traffic fatalities with other (potentially vulnerable) traffic participants by using relatively small lateral deviations to realize collision-free behavior. In this work, a complete and tractable software architecture is presented for such an AES system, comprising of the estimation of the vehicles capabilities, planning a set of paths which exploit these capabilities, checking the feasibility and risk of these paths and eventually triggering the decision to drive along one of these paths, i.e., initiating an AES manoeuvre. A novel methodology is provided to trigger such an AES system, which is based on a time-to-evade (TTE) notion. In the presence of time-varying uncertainties or measurement inaccuracies, the system is able to replan the path from the previously chosen path to ensure collision-free behavior. The proposed architecture and control approach is validated using a simulation study and field tests, showing the effectiveness of the architecture and its sub-components.  
\end{abstract}

\section{Introduction}
Autonomous Evasive Steering (AES) is a new Advanced Driver Assistance System (ADAS) feature that enhances the vehicle safety, along with existing systems like Automatic Emergency Braking (AEB), Lane Keeping Assist (LKA), Lane Departure Warning (LDW), Forward Collision Warning (FCW) and Emergency Stop Assist (ESA). This system is recommended by Euro-NCAP for new cars. Like AEB, the primary function of AES is to detect imminent collisions with other vehicles, cyclists or pedestrians based on advanced sensors and fusion of their measurements (e.g. radars, cameras and sonars). When the AES system detects an imminent collision it controls the vehicle laterally along an evasive path using the steering system and/or brake system. This system has been studied in both academic and industrial projects. One can account for the very first projects like Proreta (by Darmstadt University of Technology and Continental Automotive Systems \cite{Isermann:2008}). Since then, research and development was done in many R\&D programs by several Tier-1 and OEMs like Bosch (2010), Nissan (2012), Continental (2012) and ZF (2017) (see \cite{Dang2012} and references therein). In 2017, the first commercialized AES system was introduced on a Mercedes E-Class and on a BMW 7-series with driver-initiated versions. Lately, AES systems can also be found in the newest Ford models like the 2019 Edge, 2020 Explorer SUV \& Crossovers and some Toyota/Lexus vehicles.

\section{Previous work}
The AES-related systems have been partially or fully designed in numerous works in literature, from perception and/or risk assessment to vehicle and actuator control. Here-after, the most relevant works will be highlighted. \\
For risk assessment related to collision between vehicles and vulnerable participants (like pedestrian, cyclist), there has been several studies. For example, energy impact minimization of collision between ego vehicle and a leading car or a tailing striking car was proposed in \cite{Giugliano2015}. In \cite{Gallen2013}, collision risk is classified in distinct classes: light, serious and fatal, using the relationship between speed and car driving injury. In \cite{Naoyuki:17}, a risk map indicating a range within which the obstacle can exist after one unit of time and degree of risk based on the speed and direction of both ego and the obstacle. This risk map is then used to inform the driver or to control the vehicle in an emergency situation. In~\cite{Lee2019}, a model-predictive control (MPC) approach is used to compute the severity related to the potential crash. In \cite{serafimguardini:22}, a risk estimation based on injury curves and probabilistic occupation grid is proposed to better contextualize a scene contaning different types of objects (e.g. pedestrians, cars, cyclists ect.).\\
For path planning, many approaches have been proposed for autonomous driving, with possibility to extend to collision avoidance features. The most important examples of these approaches are stochastic methods (like RRT, MMPI \cite{Williams:2016}, SST~\cite{Smit:2022}, graph search based methods like \cite{Madas:2013}) or deterministic ones with a cost function to optimize (like MPC \cite{Funke:17}, \cite{Ploeg:2022}). Among these methods, clothoid-based path-planning is quite often used for automotive systems because of its simplicity and compatibility to human driving. In \cite{Sven:17}, clothoid curves are constructed from a known initial and final state (defined by the vehicle's states) to find a simple way to navigate from one point to another. In \cite{Lima:2015}, an MPC approach is used to optimize clothoid-based trajectories for autonomous driving.\\
For the control problem, i.e., controlling the vehicle along a trajectory, many control approaches have been proposed for collision avoidance systems. One can mention the open-loop control (\cite{Hayashi:2012}, \cite{Brannstrom:14}) up to superior closed-loop options (like PID \cite{Ackermann:14}, \cite{Sim:17}, \cite{Arbitmann:13}) or $H_\infty$ static state-feedback \cite{Hong:2017}, robust $LPV$ control \cite{DO:2021} or more modern control approaches like MPC  \cite{Funke:17}, \cite{Berntorp:17}).

Although the papers mentioned above address the concept of a collision avoidance system, they focus only on one or two main functionalities of the system. In this paper, we propose a full functional architecture and implementation of an autonomous emergency steering system. There are some works related to our proposal but to the best of our knowledge, the number of this kind of publications is rather limited. The works in~\cite{Keller:2011} and \cite{Dang2012} are amongst the first papers to propose a full concept of collision avoidance system. The designed systems in these works combine sensing, situational analysis, decision making and vehicle control (using both automatic braking and steering). However, the proposed systems are designed only for active pedestrian safety. Moreover, the path-planning algorithm, based on polynomial models, is quite simple and not a focus of the paper. The details
of the main functionalities are not provided either.

In this paper, the chosen architecture and sub-components for collision avoidance by steering is inspired from the proposed architectures for automated vehicles in \cite{Broggi:2013}, \cite{Ziegler:2014}, \cite{Gonzalez:2016} and references therein. The main contributions of the paper are the following.\\
{
\textbf{Contributions}:
\begin{enumerate}
    \item A computationally tractable and complete architecture which entails solving the path planning and control problem for an AES system while taking into account vehicle constraints
    \item A novel (re)planning concept for clothoid-based path planning
    \item A novel triggering concept for executing the AES manoeuvre
\end{enumerate}}
The outline of this article is as follows. In this section we outlined a short overview of the state of the art for AES  and collision avoidance methods. Then, we introduce the usecases and the proposed architecture in section III. In section IV, we provide technical details on all subcomponents of our proposed architecture. In section V, we back our proposed approach with simulation and experimental results and finally section VI concludes the work.
\section{Problem definition and proposed architecture}
\begin{figure}[t]
    \centering
    \includegraphics[width=\columnwidth]{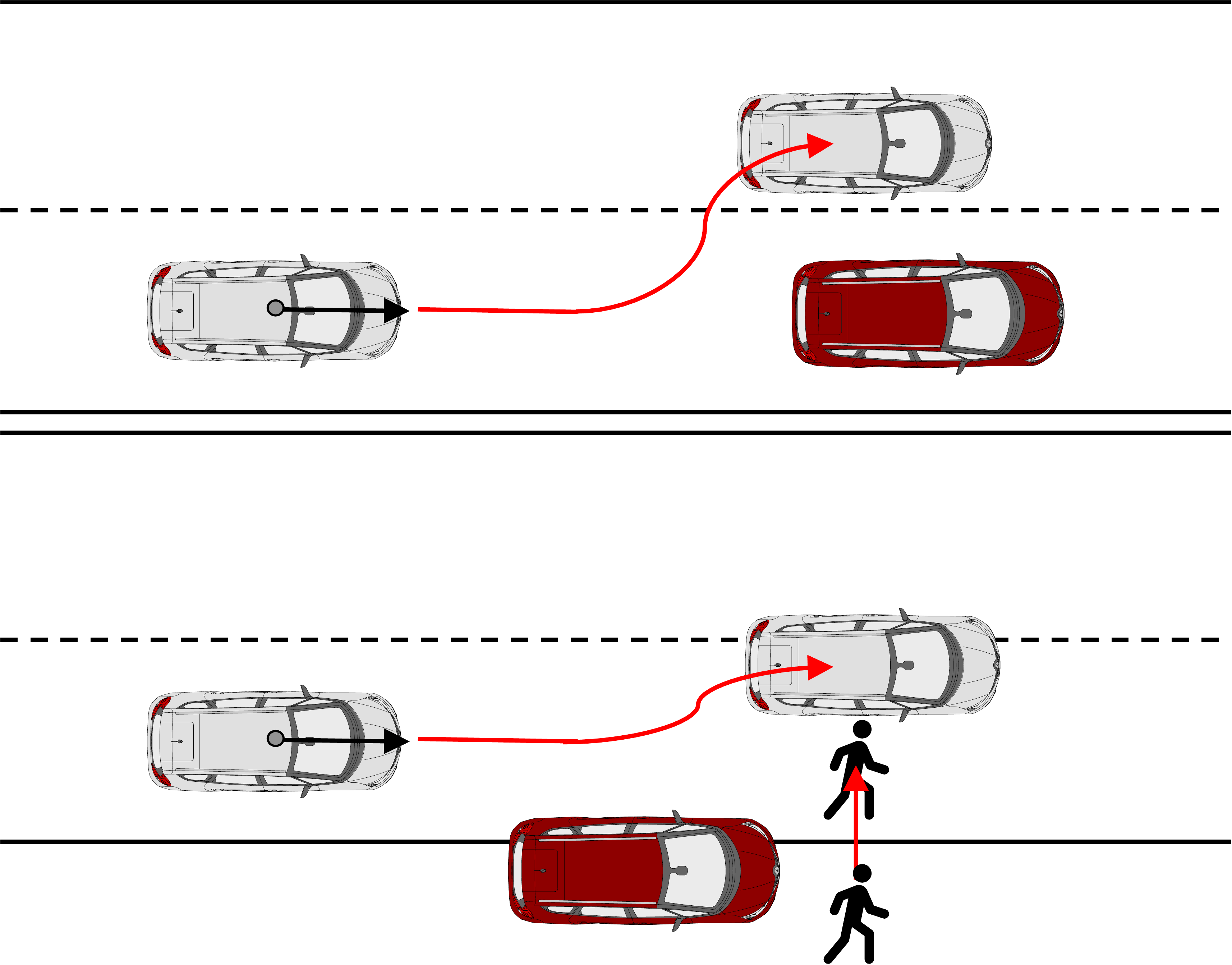}
    \caption{Use-cases in which an AES system can mitigate dangerous scenarios where an AEB system may fail.}
    \label{fig:UCselection}
\end{figure}

\begin{figure}[t]
\fontsize{8pt}{8pt}\selectfont
\def\svgwidth{\columnwidth}
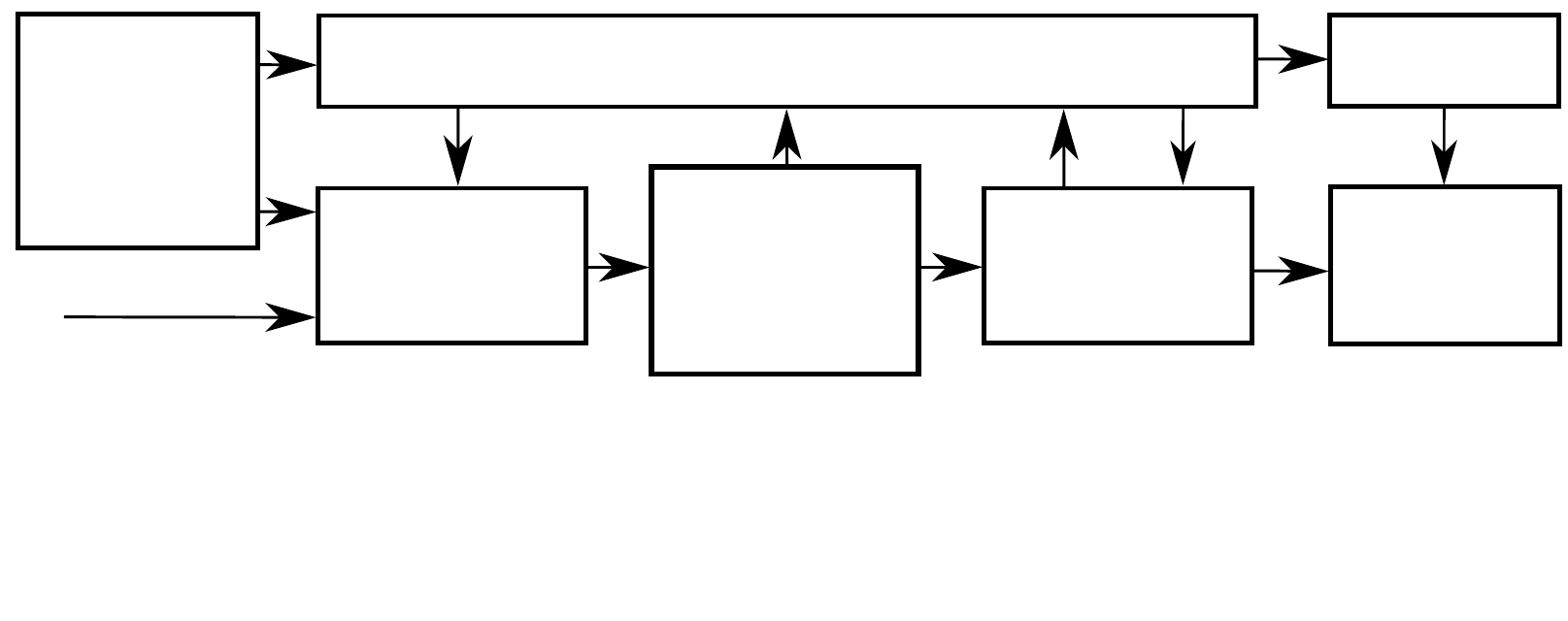
    \caption{Proposed architecture for the Autonomous Evasive Steering system.}
    \label{fig:architecture}
\end{figure}

The work in~\cite{Dang2012} motivates the use of an AES system by the reduced distance such a system needs to evade an object compared to a conventional AEB system. It states that the potential of an AES system primarily lies in mitigating accident scenarios with standstill or slowly moving objects, for which the look-ahead distance of the on-board perception is shorter than the distance needed to safely perform an AEB manoeuvre. A different category of accidents, for which an AES system is beneficial, is the case where a slow moving Vulnerable Road User (VRU), e.g., a pedestrian or a cyclist, crosses the road. In these accidents, an AEB system is only able to operate correctly when it is able to perceive the VRU from a sufficiently long distance (i.e., the braking distance needed). In case of a temporary obstruction of the VRU, i.e, the VRU appears from behind a standstill object, or the sensors are partially occluded (e.g., during bad weather) the AEB system may no longer be able to mitigate all risk.  These accidents form the basis of our research and are depicted in Figure~\ref{fig:UCselection}. 
The proposed architecture for an autonomous emergency steering system is depicted in Figure~\ref{fig:architecture}. This architecture ensures the availability of a minimal-risk collision-free trajectory at any time while allowing replanning of the trajectory in case the current triggered trajectory is no longer feasible. The main building blocks and preliminaries/assumptions are described in this section. After that, each building block is discussed in more detail in the next section.
The main starting point for the generation of dynamically feasible trajectories is to determine the current vehicle dynamic capabilities, i.e., \textit{Capability estimation}. This functionality uses measurements and/or estimates of the vehicle states (e.g., yaw rate, friction coefficients) to determine the maximal allowed lateral and longitudinal acceleration throughout the autonomous emergency steering maneuvre. The vehicle capabilities are subsequently given to the \textit{Path generation} functionality, which, by design, provides a set of dynamically feasible trajectories which are scaled within the width of the driveable space, assumed to be measurable and given. The set of paths is provided to the \textit{Path rejection/ranking} functionality, which tests whether the trajectories are enclosed by the driveable space and whether they are collision-free from dynamical objects at the current time up to the length of the trajectories. As such, we assume that the measurement and prediction of dynamic objects are available. Furthermore, the functionality ranks the driveable and collision-free paths based on severity (e.g., maximum accelerations) and proximity to the object, such that an appropriate trade-off between safety and comfort can be made. The collision-free paths are provided to the \textit{Path selection} block, which decides on the most appropriate trajectory, given the ranking appointed to them. The outcome is a collision-free minimal-cost trajectory that can be initiated by the \textit{Triggering} function which determines the last possibility to evade to prevent a collision. Once initiated the \textit{Motion control} starts controlling the vehicle along the desired trajectory through combined steering and differential braking. The \textit{State control} block controls the internal state and state transitions of the AES system. 

Practically, AES should interact with other active safety systems (e.g., Autonomous Emergency Braking (AEB), Electronic Stability Program (ESP), Anti-Lock Braking System (ABS)) to complement one another, or nominal systems (e.g., lane-keeping) to make use of them to ensure a comfortable and safe transition between the interrupts. The interactions between these systems, the vehicle itself, and the driver are considered relevant, yet out of scope for this work.
\section{Technical details}
\subsection{Capability estimation}
\begin{figure}[t]
\includegraphics[width=\columnwidth]{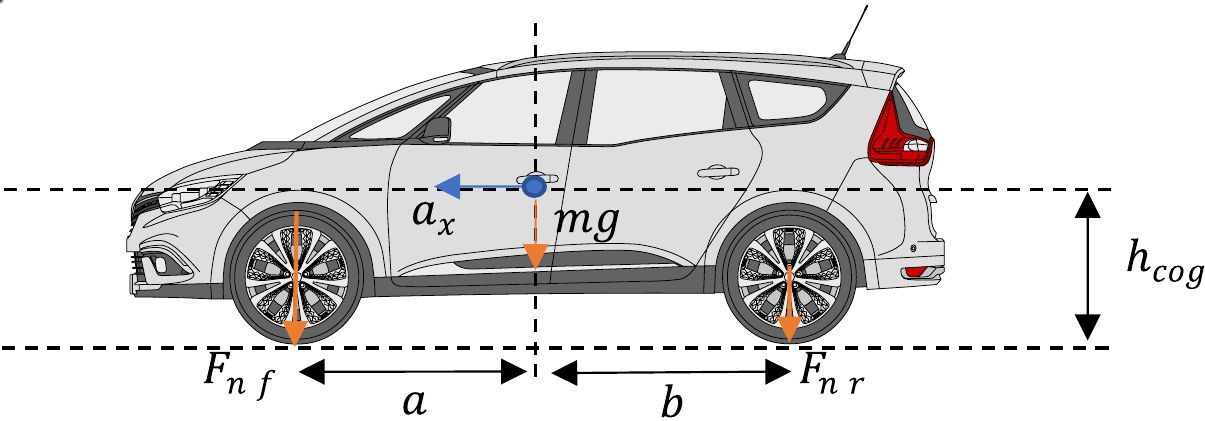}
\caption{Longitudinal free-body diagram}
\label{fig:longFBD}
\end{figure}

The \textit{capability estimation} block determines the manoeuvring capabilities of the vehicle, based on its current vehicle state and the states of the vehicle’s systems. It describes the dynamic capability along two axis: longitudinal and lateral. The longitudinal capability is used to estimate the maximum braking deceleration which can be used before starting to steer during the evasive manoeuvre. The lateral capability is used to determine the extrema of the clothoid shape of the evasive paths, i.e., the maximum steady-state curvature $\rho_{max}$ and the maximum curvature-rate $\dot{\rho}_{max}$. The lateral motion of the vehicle is performed by steering, differential braking, or a hybrid combination of the two. The capabilities are calculated for six different scenarios as depicted in Table~\ref{table:CapabilityEstimation}.
\begin{table}
\centering
\begin{tabular}{ |c|c|c|c| } 
    \hline
    Scenario & Long. braking & Steering & Diff. braking \\
    \hline
    1 & \checkmark & \checkmark & \text{\sffamily X} \\ 
    2 & \checkmark & \text{\sffamily X} & \checkmark \\ 
    3 & \checkmark & \checkmark & \checkmark \\ 
    4 & \text{\sffamily X} & \checkmark & \text{\sffamily X} \\ 
    5 & \text{\sffamily X} & \text{\sffamily X} & \checkmark \\ 
    6 & \text{\sffamily X} &\checkmark & \checkmark \\ 
    \hline
\end{tabular}\\
\caption{Different scenarios to calculate the vehicle capability.}
\label{table:CapabilityEstimation}
\end{table}
The maximum longitudinal deceleration is calculated using Newton's second law of motion, resulting in
\begin{align}
    a_{x,min}=\frac{F_{x,min}}{m},\label{eq:capability0}
\end{align}
where $m$ represents the known vehicle mass and $F_{x,min}$ represents the maximum braking force which can be calculated as follows
\begin{align*}
    F_{x,min}=-\mu_{f}F_{n,f}S_f-\mu_{r}F_{n,r}S_r,
\end{align*}
where $\mu_f,\mu_r$ represents the minimum friction coefficient of the front wheels and rear wheels, respectively, which is assumed known. Furthermore, $F_{n,f}, F_{n,r}$ represents the normal force on the front and rear axle respectively. Finally, $S_f,S_r$ denotes the effectiveness value of the braking system of the front and rear axles. For example, in case of a failure of the front braking system, $S_f=0$, in case of a loss of effectiveness, it can take any value between $0$ and $1$. Finally, the front and rear axle normal forces are calculated as follows
\begin{align*}
    F_{n,f}=&\frac{a}{a+b}mg-\frac{h_{cog}}{a+b}ma_x,\\
    F_{n,r}=&\frac{b}{a+b}mg+\frac{h_{cog}}{a+b}ma_x,
\end{align*}
where $g$ represents the gravitational constant, $a_x$ the current measured longitudinal acceleration, $a,b$ the front and rear axle distance to the vehicle center of gravity and $h_{cog}$ the center of gravity height. The involved variables are shown in Figure~\ref{fig:longFBD}. 
The lateral steering capability is determined using the maximum steady-state curvature which can be achieved using the following steady-state cornering relation obtained from the dynamic bicycle model~\cite{Jonasson2017}
\begin{align}
    \rho_{max,steering}=&\frac{\lvert\delta_{max}\rvert}{l+\frac{K_{us}v_x^2}{g}},\label{eq:capability1}\\
    K_{us}=&\frac{m}{l}\left(\frac{b}{C_f}+\frac{a}{C_r}\right),\nonumber
\end{align}
where $v_x$ represents the longitudinal velocity, $l$ represents the wheelbase, i.e., $l=a+b$. The constants $C_f,C_r$ represent the front and rear cornering stiffness. Finally, $\delta_{max}$ represents the maximum steering angle, which can be considered a tuning variable and could be obtained by, e.g., a controllability study or safety analysis (e.g., Hazard Analysis and Risk Assessment (HARA)). The lateral capability through differential braking is determined through the following steady-state relationship~\cite{Jonasson2017}
\begin{align}
    \rho_{max,diff}=\frac{w(C_f+C_r)\min{(\mu_f,\mu_r)}mg}{4\left(C_fC_r-mv_x^2\left(bC_r-aC_f\right)\right)},\label{eq:capability2}
\end{align}
where the variable $w$ denotes the track width of the vehicle. The capability for combined steering and differential braking is simply denoted as the sum of~\eqref{eq:capability1} and~\eqref{eq:capability2}. Note, that for all lateral capabilities, the current longitudinal velocity $v_x$ is used as an input variable. As mentioned before and denoted in Table~\ref{table:CapabilityEstimation}, the capabilities are also analyzed with a short braking interval (from here on denoted as "pre-braking") before the evasive manoeuvre. The pre-braking manoeuvre is introduced to induce additional load-transfer towards the front axle to reduce any potential effect of understeer due to an underestimation of the friction coefficient. The capability of this braking interval is calculated in~\eqref{eq:capability0}. Assuming this capability is used, one can directly derive the velocity $v_x$ to be used to calculate the lateral capability, as follows
\begin{align}
    v_x=a_{x,min}t_b+v_{x,0},\label{eq:prebraking}
\end{align}
where $v_{x,0}$ denotes the initial velocity, and $t_b$ denotes the duration of braking, which is considered a tuning variable. 

The lateral capabilities~\eqref{eq:capability1}, \eqref{eq:capability2} do not incorporate the fact that a certain curvature may not be achieved due to a lack of normal force and/or friction of the tyres, as such, it should be saturated through the physically maximum curvature, calculated as follows
\begin{align}
    \rho_{max,fric}=\frac{a_{y,max}}{v_x^2},\label{eq:capability3}
\end{align}
where $a_{y,max}$ denotes the maximum lateral acceleration, and is calculated as
\begin{align*}
    a_{y,max}=\min{(\mu_f,\mu_r)}g.    
\end{align*}
Finally, to make sure that a driver is able to overrule the AES system in case of a system failure, it can be decided to limit the maximum lateral acceleration of the vehicle through a pre-defined look-up table which could be scheduled, for example, based on the velocity and certain comfort and/or controllability limits. Such a threshold, $a_{y,threshold}$, results in a maximum threshold-based curvature as follows
\begin{align}
    \rho_{max,threshold}=\frac{a_{y,threshold}}{v_x^2}.\label{eq:capability4}
\end{align}
In summary, the curvature capability while steering is given by~\eqref{eq:capability1}, for differential braking ~\eqref{eq:capability2}, for combined steering and braking one can take the sum of~\eqref{eq:capability1} and~\eqref{eq:capability2}. The capabilities can be corrected for additional longitudinal braking before the manoeuvre by including~\eqref{eq:prebraking}. Finally, the steering capability can be saturated based on a lack of friction or normal force~\eqref{eq:capability3} or a pre-defined lateral acceleration threshold~\eqref{eq:capability4}. Determination of the maximum curvature rate $\dot{\rho}_{max}$ is considered out of scope for this work, as it depends on steering and braking actuator-specific properties. Throughout this work, we assume that this quantity is constant and based on the maximum human steering rate~\cite{Breuer1998}.
\subsection{Path generation and replanning}
\begin{figure}[t]
\includegraphics[width=\columnwidth]{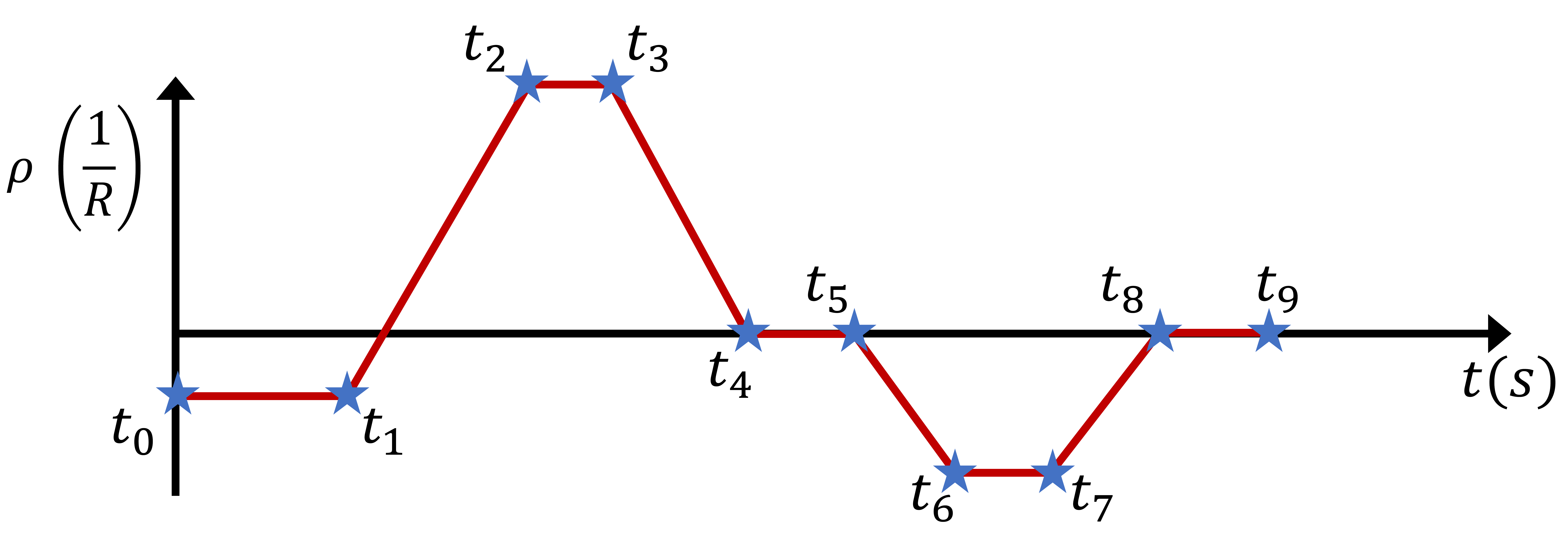}
\caption{AES path curvature model with indicated the timing of the break points with $t_0\rightarrow t_9$.}
\label{fig:PathCurvatureModel}
\end{figure}

The path generation block generates paths both for the purpose of planning an initial AES trajectory, as well as replanning throughout the duration of a triggered AES manoeuvre. An evasive action is a short high dynamic manoeuvre. This means that the planned path is required to take the dynamical behaviour and constraints of the vehicle into account in order to ensure the driveability of the path. Furthermore, an evasive action is often characterized by a distinct pattern; a lateral motion to the left or right of the obstacle to avoid a collision and a stabilizing action to put the vehicle on the original course again. As discussed in the introduction, the trajectory is planned within the action space such that dynamic feasibility of driving along the path can be ensured. To generate the evasive path, the path generator calculates up to 10 points in terms of time $t$ with the accompanying curvature $\rho$ and longitudinal velocity $v_x$. An example path is illustrated in Figure~\ref{fig:PathCurvatureModel}. Throughout the construction of action space paths,  certain tuning parameters can't be derived analytically and should be determined through qualitative arguments (e.g., safety analysis, comfort, etc.). First, the parameter $t_{pb}$ is introduced, which represents the time-interval for which a pre-braking action, as presented in~\eqref{eq:prebraking}, is sustained. Furthermore, the parameter $\psi_{max}$ is introduced, which denotes the maximum heading angle the vehicle is allowed to achieve during an AES manoeuvre. This section is divided in different subsections to explain how the data points are generated at each time sample to obtain a maximum severity AES path (i.e., satisfying all capabilities with and without pre-braking). Subsequently, the sampling of the most severe path to a set of AES paths within the driveable area is explained. Furthermore it is discussed how to incorporate road-curvature profiles to superimpose on the existing AES paths. For the sake of simplicity, it is assumed that the road-curvature is constant, an extension towards potentially variable road curvatures is straight forward.
\subsubsection{Trajectory point $t_0$}
The initial condition of the AES trajectory in action space should satisfy the initial conditions of the vehicle. That is, the curvature at $t_0$,  $\rho_0$, and the velocity at $t_0$, $v_{x,0}$, must satisfy
\begin{align*}
    \rho_0=&\frac{\dot{\psi}_{ego}}{v_{x,ego}},\\
    v_{x,0}=&v_{x,ego},
\end{align*}
where $\dot{\psi}_{ego}$ represents the measured yaw-rate of the ego vehicle and $v_{x,ego}$ represents the measured longitudinal velocity of the vehicle at current time instances.
\subsubsection{Trajectory point $t_1$}
During time interval $t_0\rightarrow t_1$, the AES manoeuvre is initiated. The curvature at $t_1$ is determined by the initial curvature $\rho_0$, corrected with the velocity as the result of a pre-braking action~\eqref{eq:prebraking}. In the absence of a pre-braking action, one can consider $a_{x,min}=0\:\text{m}\cdot \text{s}^{-2}$. The velocity, curvature and time-interval $t_1-t_0$ are hence calculated as
\begin{align*}
    t_1 = &t_{pb}+t_0\\
    v_{x,1} =& v_{x,0}+a_{x,min}(t_{1}-t_0)\\
    \rho_1 =& \frac{{v_{x,0}}\rho_0}{v_{x,1}}+\rho_{road},
\end{align*}
where $\rho_{road}$ represents the curvature of the road.
\subsubsection{Trajectory point $t_2$}
From time interval $t_1\rightarrow t_2$ the evasion is initiated. The determination of the path during this time interval is more complex, as it starts taking into account the various constraints the path generation algorithm has. From $t_1$ onwards, the velocity is considered constant, i.e., 
\begin{align*}
    v_{x,i}=v_{x,1}\quad\forall i\in[1,2 \hdots 9].
\end{align*}
First, we calculate the accumulated heading error from $t_0\rightarrow t_1$ as follows
\begin{align*}
    \psi_{tot,1}=&\int_{t_0}^{t_1}\dot{\psi}(t)dt,\\
    =&\frac{(t_1-t_0)\rho_0 v_{x,0}+(t_1-t_0)\rho_1 v_{x,1}}{2}.
\end{align*}
We assume, for now, that the initial heading ($\psi_0,\psi_1$), curvature ($\rho_0,\rho_1$) and road curvature is equal to zero. The time intervals $t_1\rightarrow t_2$ and $t_3\rightarrow t_4$ are calculated as:
\begin{align}
    \Delta t_{1\rightarrow 2}+\Delta t_{3\rightarrow 4}=& \frac{1}{2}\frac{\lvert\rho_2\rvert}{\dot{\rho}_{max}}+\frac{1}{2}\frac{\lvert\rho_3\rvert}{\dot{\rho}_{max}},\nonumber\\=&\frac{\lvert\rho_2\rvert}{\dot{\rho}_{max}}.\nonumber
\end{align}
And the time it takes, given a curvature $\rho_2$, to reach the maximum desired heading is calculated as
\begin{align*}
    \Delta t_{1\rightarrow 4} = \frac{\dot{\psi}_{max}}{\lvert\rho_2 \rvert v_{x,2}}
\end{align*}
In order to reach the maximum allowed heading $\psi_{max}$ from time-interval $t_0\rightarrow t_4$, while enforcing $\rho_2=\rho_3$, we could first try to find a $\rho_2$ such that $\Delta t_{2\rightarrow 3}=0$. Hence, we are looking to solve the following equilibrium:
\begin{align}
    \Delta t_{1\rightarrow 4} =& \Delta t_{1\rightarrow 2}+\Delta t_{3\rightarrow 4}\nonumber\\
    \frac{\psi_{max}}{\lvert\rho_2\rvert v_{x,2}}=&\frac{\lvert\rho_2\rvert}{\dot{\rho}_{max}}.\label{eq:rho2eq}
\end{align}
One can solve the above equilibrium to find a desired curvature $\rho_2$, such that the maximally desired heading is achieved. If the vehicle is initiated with non-zero heading, i.e., $\psi_{ego}\neq 0$, the equation~\eqref{eq:rho2eq} is recast as follows
\begin{align*}
    \frac{\psi_{max}-(\psi_{ego}+\psi_{tot,1})}{\lvert\rho_2\rvert v_{x,2}}=\frac{\lvert\rho_2\rvert}{\dot{\rho}_{max}},
\end{align*}
for which, again, a closed-form solution for $\rho_2$ could be obtained. Consider, the fact that $\rho_2$ represents the maximum curvature of the trajectory, we require it to not surpass the maximum capability, i.e., $\rho_{max}$. The closed-form solution for $\rho_2$ and hence $t_2$ is therefore denoted as
\begin{align*}
    \rho_2\!=\!&\min\!\left(\!\sqrt{\frac{(\psi_{max}\!-\!(\psi_{ego}\!+\!\psi_{tot,1}))\dot{\rho}_{max}}{v_{x,2}}}\!+\!\rho_{road},\lvert \rho_{max}\rvert\!\right),\\
    t_2=&\frac{\rho_2-\rho_1}{\dot{\rho}_{max}}+t_1.
\end{align*}
\subsubsection{Trajectory point $t_3$}
From time interval $t_2\rightarrow t_3$ the increase in curvature, determined in interval $t_1\rightarrow t_2$ is stopped, i.e.,
\begin{align*}
    \rho_3=\rho_2.
\end{align*}\
If, at time instance $t_2$, the curvature $\lvert\rho_2\rvert=\rho_{max}$, the time interval $t_2\rightarrow t_3$ can be utilised to reach the maximum allowable heading $\psi_{max}$, i.e.,
\begin{align*}
    v_{x,3}(t_3-t_2)\rho_3=(\psi_{max}-(\psi_{ego}+\psi_{tot,1}))-(t_2-t_1)\rho_2v_{x,2}\\
    t_3 = \frac{(\psi_{max}-(\psi_{ego}+\psi_{tot,1}))-(t_2-t_1)\rho_2v_{x,2}+t_2v_{x,3}\rho_3}{v_{x,3}\rho_3}
\end{align*}
\subsubsection{Trajectory point $t_4$}
During time interval $t_3\rightarrow t_4$, the curvature is brought back to the curvature of the road. By assuming that, for now, the road geometry is straight, the curvature is hence equal to zero, i.e., 
\begin{align*}
    \rho_4=0.
\end{align*}
The time instance $t_4$ is, hence, easily determined by
\begin{align*}
    t_4=\frac{\rho_3}{\dot{\rho}_{max}}+t_3.
\end{align*}
\subsubsection{Trajectory point $t_5$}
Time interval $t_4\rightarrow t_5$ can be used to generate a higher lateral deviation with respect to the initial state of the vehicle. This can be done, e.g., to reach the outer boundary of the driveable space. The lateral offset, generated during this time interval, can be calculated analytically by
\begin{align*}
    y_{offset}=v_{x,4}(t_5-t_4)\sin{(\psi_4)}
\end{align*}
where $t_5$ can be calculated based on the desired additional lateral offset, and the curvature $\rho_5$ is equal zo zero, i.e.,
\begin{align*}
    t_5=&\frac{y_{offset}+v_{x,4}t_4\sin{(\psi_4)}}{v_{x,4}\sin{(\psi_4)}}\\
    \rho_5=&\rho_4=0
\end{align*}
\subsubsection{Trajectory point $t_6$}
As mentioned before, the AES path phases are twofold. First, there is the evasive phase, where the vehicle moves off it's initial trajectory to evade a critical object. The evasive phase occurs from time instance $t_0$ up until time instance $t_5$. From $t_6$ onward, the objective is to stabilize the vehicle and bring the vehicle parallel to the road. In this case the heading at the final time instance should be 0, i.e., $\psi_9=0$. Hence, given that $\psi_5=\psi_4$, the following should hold
\begin{align}
    \psi_5=&\frac{1}{2}(v_{x,6}\rho_6(t_6-t_5)+v_{x,8}\rho_7(t_8-t_7))\nonumber\\&+v_{x,7}\rho_7(t_7-t_6),\label{eq:t6calc}
\end{align}
which implies that the path ends up in parallel to the road. Moreover, the following constraints hold
\begin{align*}
    \rho_7=&\rho_6\\
    t_6-t_5=&t_8-t_7
\end{align*}
furthermore, the velocity is still assumed constant and equal to $v_{x,1}$, hence equation~\eqref{eq:t6calc} is reduced to
\begin{align*}
    v_{x,6}\rho_6(t_6-t_5)+v_{x,6}\rho_6(t_7-t_6)&=\psi_5,\\
    v_{x,6}\rho_6(t_7-t_5)&=\psi_5.
\end{align*}
Where the maximum value of $\rho_2,\rho_3$ was chosen as the maximum lateral capability, the value of $\rho_6$ is defined as $\rho_6=\rho_2i_{sb}$ where $i_{sb}<1$ is a tuning parameter to prevent instability of the vehicle caused by oversteering behavior (an intuitive example of this behavior is the so-called "Scandinavian flick"). This is a concept choice, do note, that more extreme manoeuvres could be induced by choosing $i_{sb}=1$, which could require controlling the vehicle on or beyond the limit of handling. Now that an initial estimate for $\rho_6$ is known, it has to be checked whether reaching this curvature with the maximum rate of curvature is sufficient to compensate for the heading $\psi_{max}$, i.e., if
\begin{align*}
    \frac{\lvert\psi_5\rvert}{\lvert\rho_2i_{sb}\rvert v_{x,6}}-\frac{\lvert\rho_2i_{sb}\rvert}{\dot{\rho}_{max}}=0,
\end{align*}
then, the used curvature is exactly sufficient to overcome the heading without using time interval $t_6\rightarrow t_7$. However, if
\begin{align*}
    \frac{\lvert\psi_{max}\rvert}{\lvert\rho_2i_{sb}\rvert v_{x,6}}-\frac{\lvert\rho_2i_{sb}\rvert}{\dot{\rho}_{max}}
    <0,    
\end{align*}
the used curvature is too high. A sufficient curvature to overcome the heading $\psi_5$ can be calculated by solving the following equilibrium for $\rho_6$
\begin{align*}
    \frac{\lvert\psi_{max}\rvert}{\lvert\rho_6\rvert v_{x,6}}=\frac{\lvert\rho_6\rvert}{\dot{\rho}_{max}}.
\end{align*}
Finally, if the following condition holds
\begin{align}
    \frac{\lvert\psi_{max}\rvert}{\lvert\rho_2i_{sb}\rvert v_{x,6}}-\frac{\lvert\rho_2i_{sb}\rvert}{\dot{\rho}_{max}}
    >0,   \label{eq:t7condition} 
\end{align}
the used curvature is too low, meaning that an additional constant curvature phase $t_6\rightarrow t_7$ needs to be used. In summary, the curvature $\rho_6$ is calculated through the following solution
\begin{align*}
\rho_6=\min\left(\sqrt{\frac{\psi_{max}\dot{\rho}_{max}}{v_{x,6}}}+\rho_{road},i_{sb}\rho_2\right),
\end{align*}
and the time interval $t_5\rightarrow t_6$ needed to reach that curvature equals to
\begin{align*}
    t_6-t_5=&\frac{\lvert\rho_6\rvert}{\dot{\rho}_{max}}.
\end{align*}
\subsubsection{Trajectory point $t_7$}
As is mentioned in the previous section, time-interval $t_6\rightarrow t_7$ is needed if and only if the condition in~\eqref{eq:t7condition} holds. In this case, we set
\begin{align*}
    \rho_7=\rho_6,
\end{align*}
and the time instance $t_7$ is calculated as
\begin{align*}
    (t_7-t_6)=&\frac{\lvert\psi_{max}\rvert}{\lvert\rho_6\rvert v_{x,6}}-\frac{\lvert\rho_6\rvert}{\dot{\rho}_{max}},\\
    t_7=&\frac{\lvert\psi_{max}\rvert}{\lvert\rho_6\rvert v_{x,6}}-\frac{\lvert\rho_6\rvert}{\dot{\rho}_{max}}+t_6.
\end{align*}
\subsubsection{Trajectory point $t_8$}
During time interval $t_7\rightarrow t_8$, the curvature of the path transitions to the final curvature (i.e., zero on a straight road. The derivation is similar to the derivation provided for time instance $t_6$, hence
\begin{align*}
    t_8-t_7=\frac{\rho_6}{\dot{\rho}_{max}},
    t_8=\frac{\lvert\rho_6\rvert}{\dot{\rho}_{max}}+t_7,
\end{align*}
where the value of $\rho_8$ is equal to zero for a straight road.
\subsubsection{Trajectory point $t_9$}
The final time instance of the trajectory essentially forms a tuning parameter. In essence, the goal of the time interval $t_8\rightarrow t_9$ is to allow stabilization of the vehicle. After time instance $t_8$, the heading of the path is parallel to the road, hence, any time interval $t_8\rightarrow t_9$ is chosen to allow for the motion controller to converge to the AES path.
\subsubsection*{Path sampling algorithm}
\begin{figure}[t]
    \subfloat{
        \includegraphics[width=0.9\columnwidth]{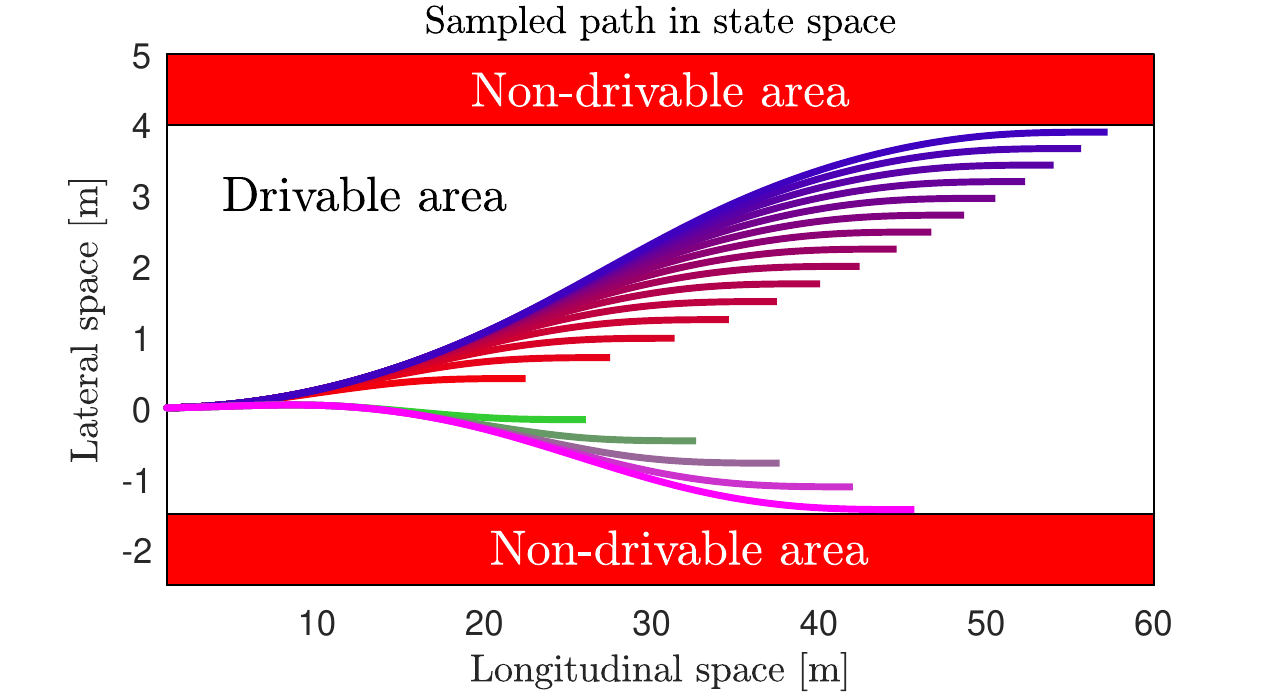}}\hfill
            \subfloat{
        \includegraphics[width=0.9\columnwidth]{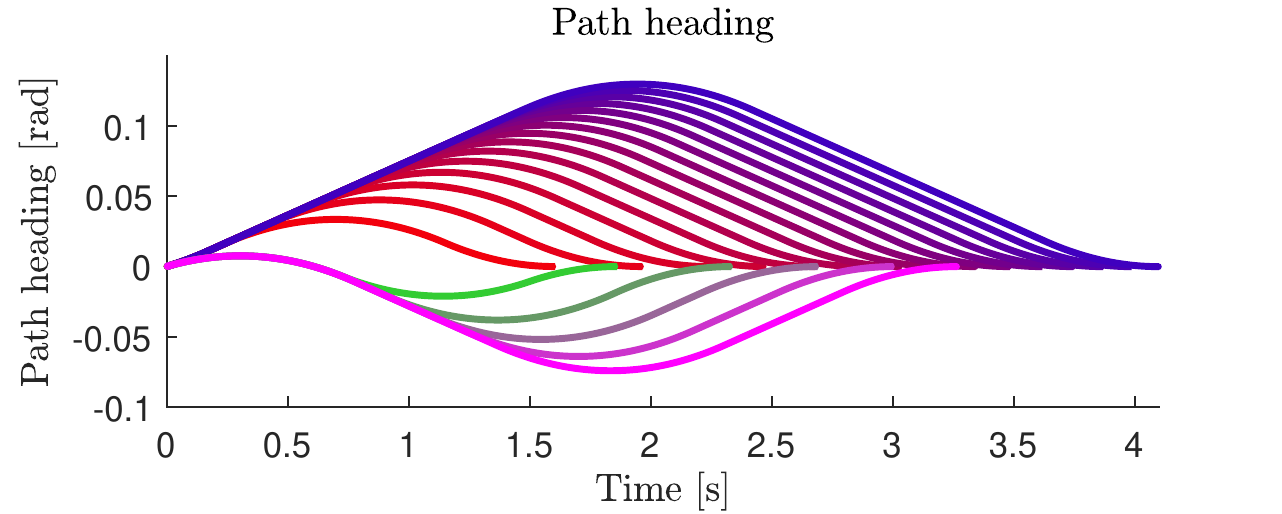}}\hfill
            \subfloat{
        \includegraphics[width=0.9\columnwidth]{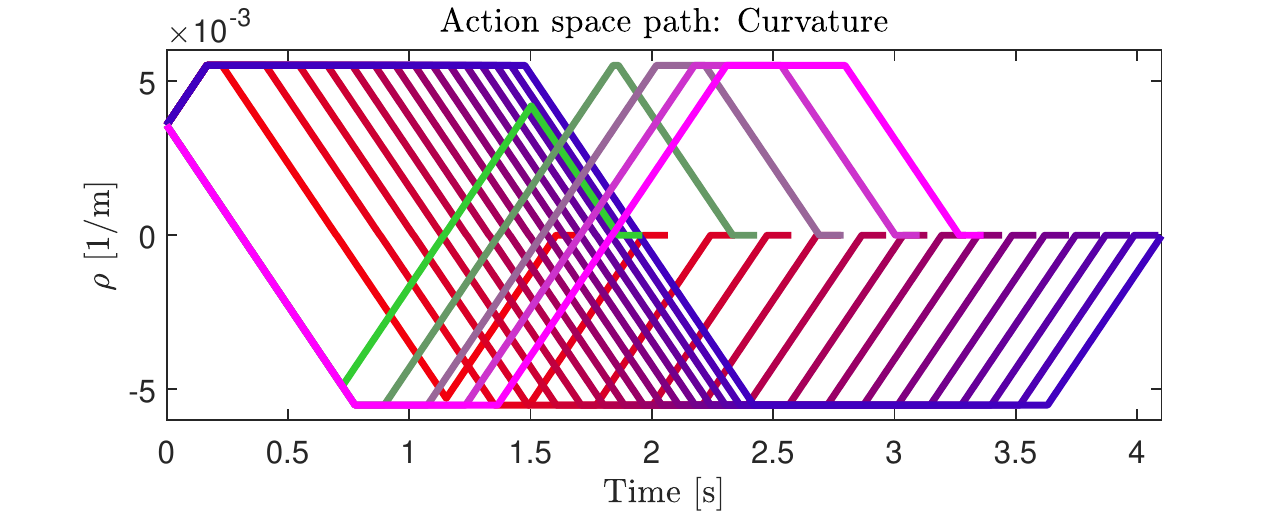}}
    \caption{Sample result of a set of generated trajectories.}
    \label{fig:pathplanningconcept}
\end{figure}

The maximum capability path, originating from the previous sections lies, by definition, on the edge of the driveable area or beyond the driveable area. As a result, it is desired to scale back this maximum capability path to fit within the driveable area. The lateral deviation achieved by travelling across a path can be found using the discrete Fresnel integral
\begin{align*}
    \psi(k)=&\psi(k-1)+\rho(k)v_x(k)\Delta t,\\
    y(k)=&y(k-1)+\sin(\psi(k))v_x(k)\Delta t.
\end{align*}
When using a sufficiently small timestep $\Delta t$, an accurate lateral deviation can be derived of the maximum capability path. We call this processing step "pre-sampling". If the lateral deviation of the path goes beyond the driveable area, the capabilities of the desired set of paths is defined through the following heuristically found relationship
\begin{align}
    \rho_{max,n}=\frac{y_{desired}}{y_{max}}\rho_{max}\sqrt{\frac{n}{n_{tot}}}\forall n\in[1,2,\hdots,n_{tot}],\label{eq:sampling1}\\
    \psi_{max,n}=\frac{y_{desired}}{y_{max}}\psi_{max}\sqrt{\frac{n}{n_{tot}}}\forall n\in[1,2,\hdots,n_{tot}].\label{eq:sampling2}
\end{align}
where $n_{tot}$ represents the desired number of paths. By using the relationships~\eqref{eq:sampling1}~\eqref{eq:sampling2}, an array of paths can be found which covers the left-hand side or the right-hand side of the vehicle with dynamically feasible paths, fitting in the left-hand or right-hand driveable space, respectively. A sample result is provided in Figure~\ref{fig:pathplanningconcept}, where all trajectories respect the driveable area, while the heading is saturated by a pre-defined maximum, and the curvature and its derivative is limited.

\subsubsection*{Replanning concept}
\begin{figure}[t]
    \subfloat{
        \includegraphics[width=0.9\columnwidth]{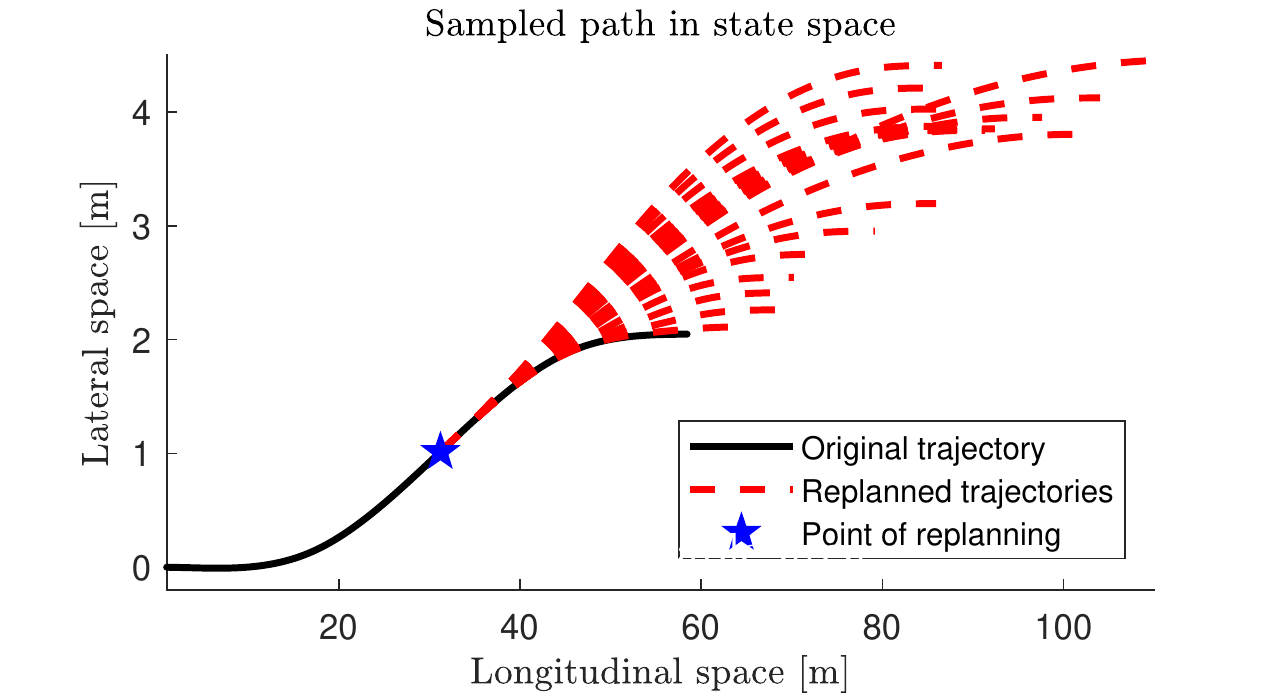}}\hfill
            \subfloat{
        \includegraphics[width=0.9\columnwidth]{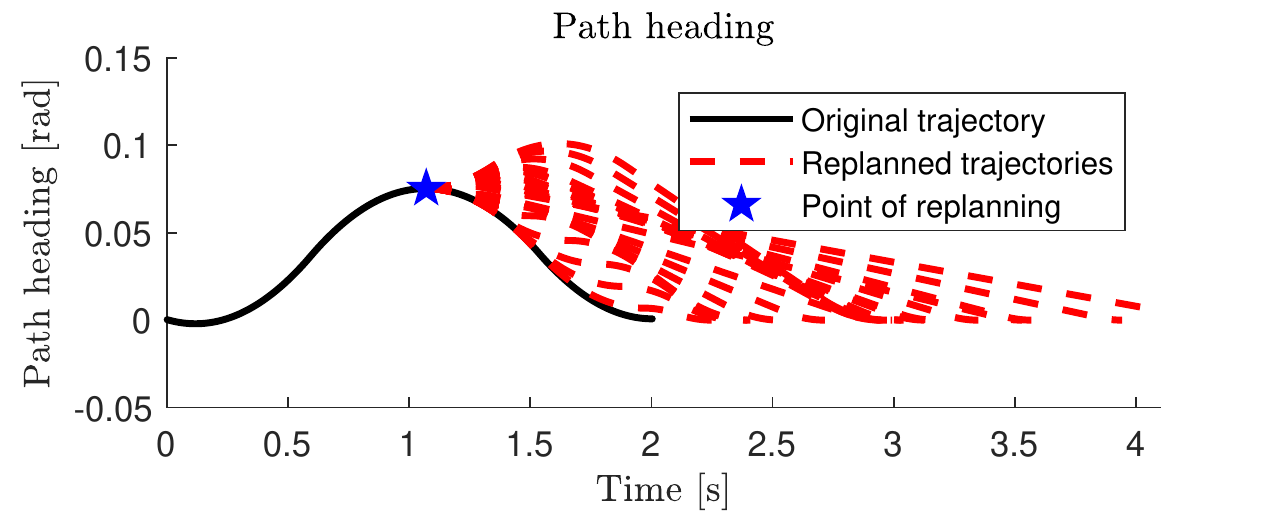}}\hfill
            \subfloat{
        \includegraphics[width=0.9\columnwidth]{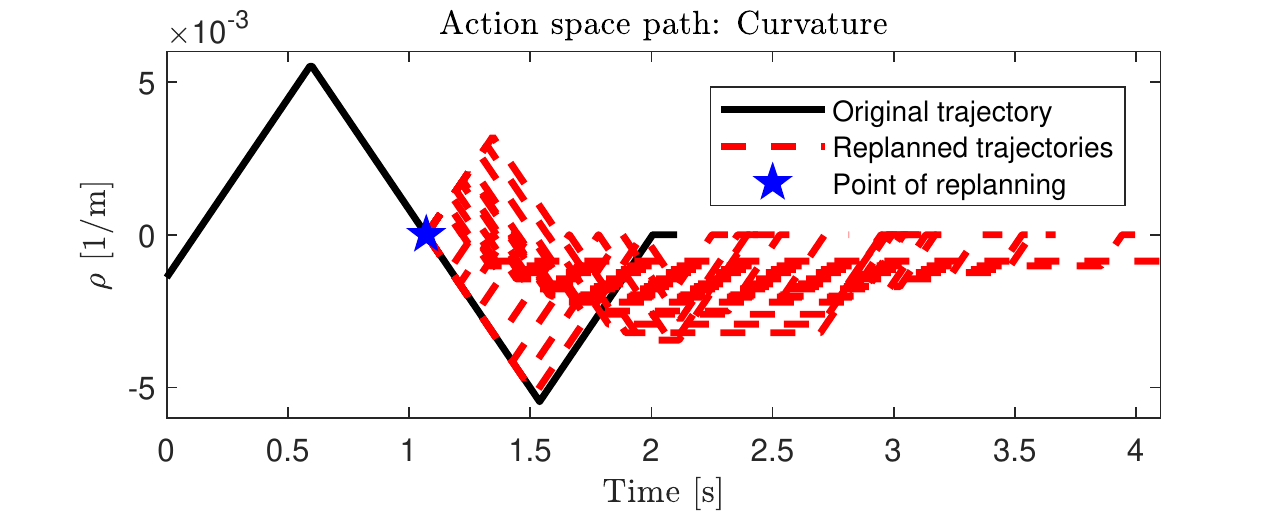}}
    \caption{Sample result of a set of replanned trajectories.}
    \label{fig:replanningconcept}
\end{figure}

In the previous sections, a methodology has been described which plans trajectories from an initial state to a state with a lateral offset within the driveable area, while making use of, and respecting, the vehicle dynamic capabilites of the vehicle. At any time while traversing an AES trajectory, it can occur that the trajectory is no longer feasible (e.g., due to new, previously unseen, observations from the vehicle perception). The requirements of a replanned path are identical to the requirements of a "regular" path before triggering the AES manoeuvre. As such, the complete calculation method for generating a set of dynamically feasible paths from $t_0\rightarrow t_9$ can be re-used by updating the initial condition at $t_0$ with the current vehicle state, and updating the driveable space. Results of the generation of re-planned paths from a chosen moment in time can be observed in Figure~\ref{fig:replanningconcept}
\subsection{Path rejection and ranking}
\begin{figure}[t]
\includegraphics[width=\columnwidth]{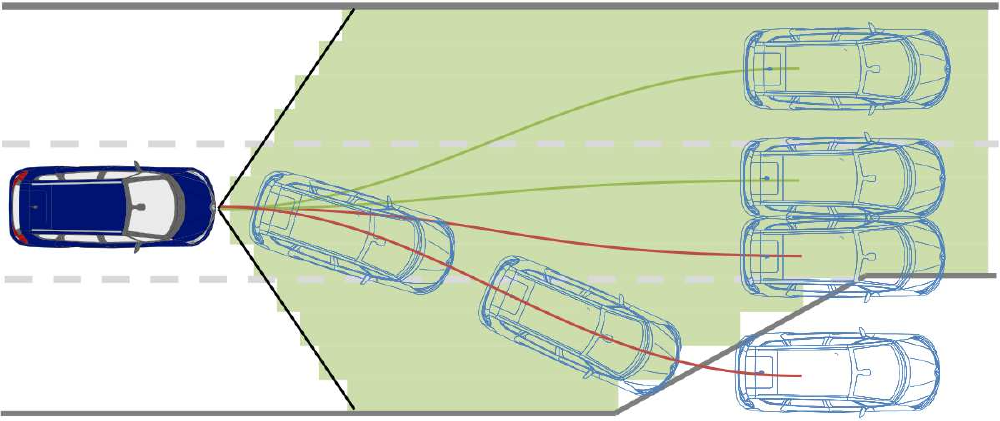}
\caption{Example of the driveable space modelling (green horizontal bars) in a specific scenario of a merging lane.}
\label{fig:DriveableSpace}
\end{figure}

\begin{figure}[t]
\includegraphics[width=\columnwidth]{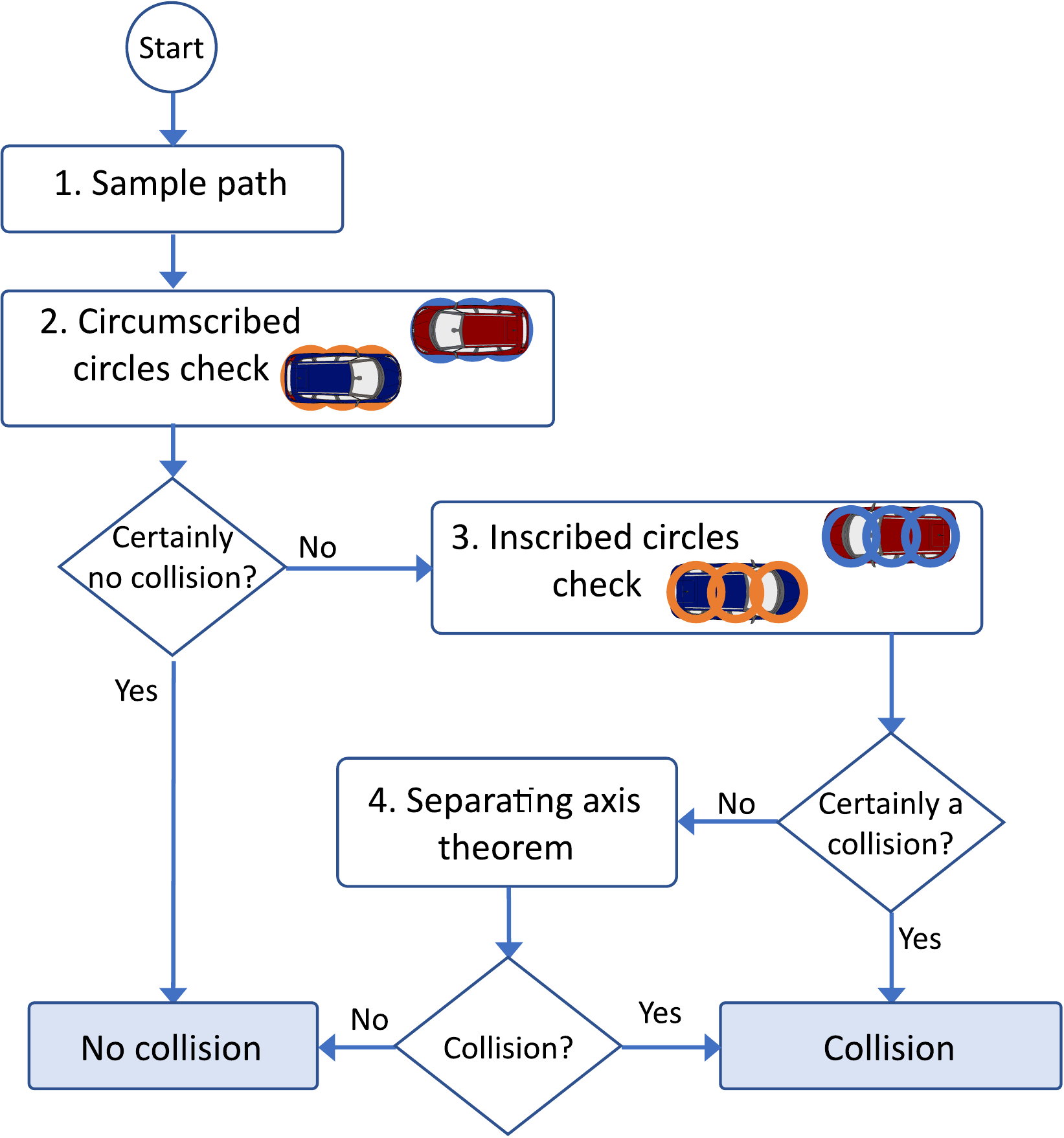}
\caption{Four step approach to check the path for collisions.}
\label{fig:StepApproachCollisionCheck}
\end{figure}

The two main objectives of the path rejection and ranking module are to reject paths that are not collision-free and calculate the cost for the remaining subset of paths. The following approach is used. First, it is checked whether the path lies within the \underline{driveable area}. Subsequently a \underline{check on collisions} with other occupants is performed. For all paths that are not rejected because they fail the driveable area and/or collision check, a \underline{cost is calculated} based on path severity and proximity to targets.
To perform the checks and to allow the motion controller to follow the path, the path received from the path generator needs to be sampled. First all paths are sampled to a rough grid to perform the checks. The selected path is sampled to a fine grid at the end, for the motion controller to use. Next to the checking of the new paths, also the path used during the manoeuvre is repeatedly checked for validity.
Checking whether the paths are within the driveable area is done by a geometric check that verifies that the swept path of the vehicle along the paths is within the driveable space that is defined along the road geometry. To do this efficiently, in the lateral direction the driveable space is divided into segments relative to the vehicle of approximately 0,5 meter wide (see~Figure~\ref{fig:DriveableSpace}.).
To limit the calculation effort of the collision check it is done according the four stpdf depicted in~Figure~\ref{fig:StepApproachCollisionCheck}. 
The collision check requires prediction of both the ego vehicle motion along the path and the targets. For this, all of the generated path are sampled (1) in time. The number of samples is a tradeoff between accuracy and calculational effort. At each time instance, first, a circumscribed circle check~\cite{Stiller2010} (2) (i.e., checking if the circles, covering the exterior of the ego vehicle and the target overlap) is performed to identify instances at which definitely no collision will occur. For instances for which collision cannot be excluded an inscribed circle check~\cite{Stiller2010} (3) (i.e., checking if the circles, covering the interior of the ego vehicle and target contour overlap) is done to identify instances at which definitely a collision will occur. For the remainder (doubtful) instances the separating axis theorem (4) is used to determine whether a collision is inevitable.

To select one of the remaining paths, a cost for every path is calculated. This cost incorporates a path severity and proximity to targets component. 
\begin{align}
    C_{path}=C_{severity} + C_{proximity},
\end{align}
Since the description of the path only gives the time, curvature, velocity and positions, these are calculated as follows:
\begin{align}
    C_{severity}=& K_{ay}\sqrt{\sum_{n=1}^{N} \lvert{v^2_{x,path}(n) \cdot \rho_{path}(n)\lvert^2}}\\& +\nonumber K_{ax}\sqrt{\sum_{n=2}^{N} \lvert \frac{v_{x,path}(n)-v_{x,path}(n-1)}{t(n)-t(n-1)}\lvert^2}
\end{align}

\begin{align}
C_{proximity}\!\!=\!\!K_{prox} \frac{\sum_{n=1}^{N} \!\min(d_{tgt}(n,i))}{N}\:\forall i \in\! (1,N_{tgt})
\end{align}
Where $N$ is the number of sampled path points, $v_{x,path}$ the velocity along the path and $d_{tgt} (n,i)$ is the distance to target $i$ at prediction step $n$. The weighing factors that allow tuning of the individual parameters are $K_{ay}$ , $K_{ax}$ and $K_{prox}$. Similar to a model-predictive formulation~\cite{Ploeg:2022}, these weights emphasize the importance of certain design-criteria for an "optimal" AES path. For example, in the presence of uncertainties, it may be preferable to increase the weight or importance of the proximity, $K_{prox}$, at the cost of a higher path severity in longitudinal and lateral sense. Moreover, it may be preferable to penalize a high lateral acceleration through $K_{ay}$, with respect to the longitudinal acceleration, in order to improve the controllability of the vehicle if a driver decides to intervene. This modular construction of weighting states of interest also allows introduction of other cost functions such as use of less preferred driveable area (i.e., generating a path which drives over a curb, but kepdf the vehicle and its surroiundings safe).

\subsection{Path selection}
Path selection function exists of three main functions; select the lowest cost path, sample the path and monitor the validity of the previous path. First the lowest cost path is selected from the collision free path. Subsequent this path is resampled in order to provide a more detailed path to the motion control. Finally the last selected path (maybe already initiated) is checked on collisions to facilitate the decision functions to replan if it results in a collision due to unforeseen change of the situation.  

\subsection{Evasive trigger}
\begin{figure}[t]
\includegraphics[width=\columnwidth]{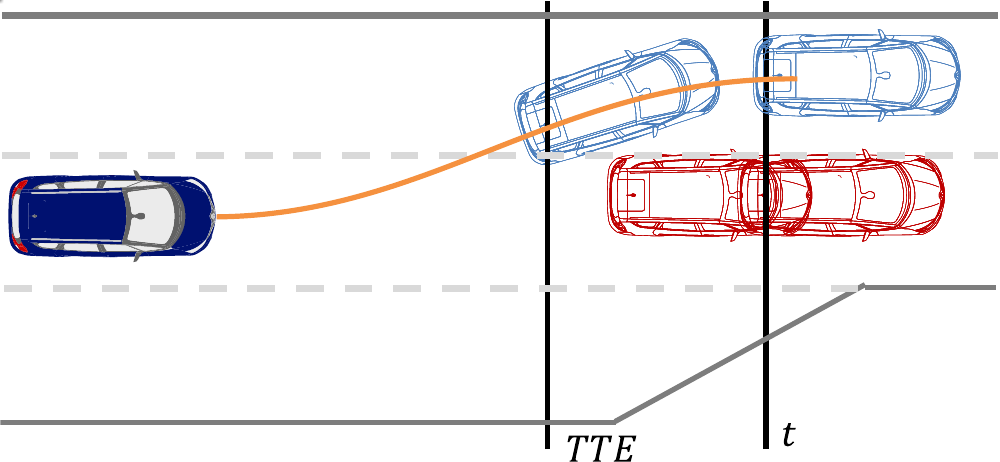}
\caption{Time To Evade (TTE) is determined by reducing t8 with a constant tunable factor.}
\label{fig:TTE}
\end{figure}

The evasive trigger, i.e., the condition to transition to an AES manoeuvre, is based on the last possibility to evade and prevent a collision. Other solutions for the trigger could be based on the Time-to-collision (TTC) or on the cost of a specific path. A trigger based on a TTC is often used for systems like AEB or collision warning systems. These systems, however, are one dimensional since they only take longitudinal braking into account. For AES, also lateral motion needs to be taken into account, and thus all relevant objects and the driveable area should be incorporated in the decision. This implies that the trigger should be based on the scenario, instead of the TTC to one specific target of interest. Basing the trigger on a cost function only is not very insightful since the cost is a combination of various factors with different magnitudes and units.
The path generation curvature is shaped according the model shown in~Figure~\ref{fig:PathCurvatureModel}. The advantage of using this model is that the time at which the evasive manoeuvre is completed ($t_8$) and a stabilization phase ($t_8\rightarrow t_9$) is initiated is known and available. To compensate for the fact that the collision may already be mitigated before the manoeuvre is completed, $t_8$ is reduced with a tunable factor to calculate the Time To Evade (TTE). This is illustrated in~Figure~\ref{fig:TTE}. 
Using the TTE and the TTC the AES triggering function will provide a warning to the vehicle and driver when the following inqeuality holds:
\begin{align}
    TTC \leq TTE + t_{margin}+ t_{warning}
    \label{eq:Warningtrigger}
\end{align}
where $t_{warning}$ is the warning time, and $t_{margin}$ is the trigger window (i.e., the margin on top of the TTE to prevent triggering the AES too late). Note, that the tuning variable $t_{margin}$ can also be used to compensate for potential actuation and/or sensing delays. When the TTC decreases further and the following inequality holds, AES will transition into the AES In Regulation state and will initiate the evasive manoeuvre: 
\begin{align}
    TTE \leq TTC \leq TTE + t_{margin}
    \label{eq:TimeToEvade}
\end{align}

\subsection{State control}
\begin{figure}[t]
\fontsize{8pt}{8pt}\selectfont
\def\svgwidth{\columnwidth}
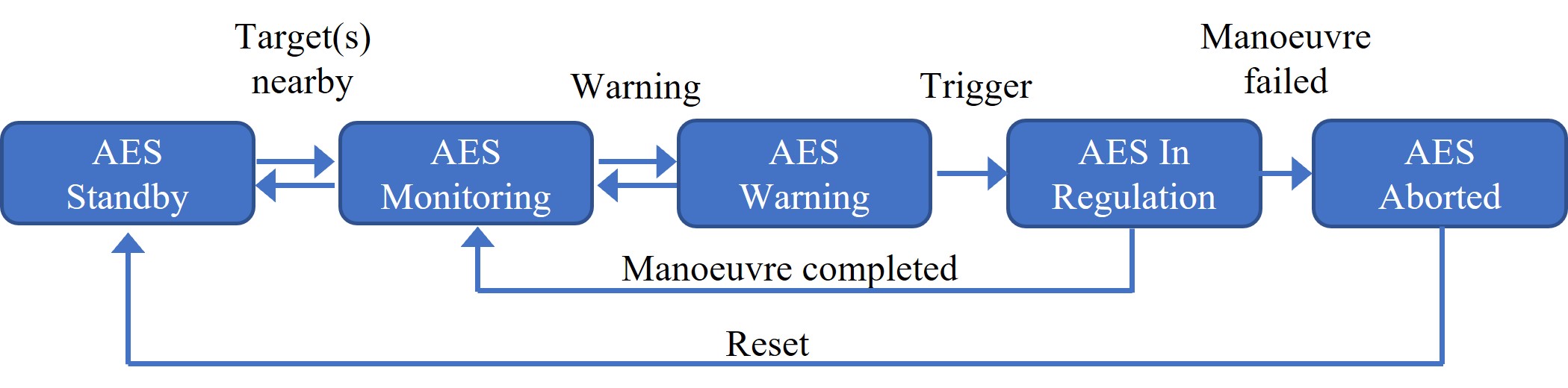
    \caption{AES state diagram including the transition conditions.}
    \label{fig:AESStateDiagraml}
\end{figure}

To orchestrate the evasive manoeuvre a state controller is added. In~Figure~\ref{fig:AESStateDiagraml} the basic state flow diagram of this controller is provided.
When initialised, the AES system will go into the \textit{AES standby} state. When target information is made available by the perception function the state will change to \textit{AES monitoring} in which planning of evasive paths is enabled. When all targets are out of range or the perception function fails, AES will go back again to the \textit{AES standby} state. When the monitoring state equation~\eqref{eq:Warningtrigger} is satisfied and there exist no system errors preventing execution of an evasive manoeuvre then an AES warning to the driver will be issued. If the TTC increases, then the state can change back to \textit{AES monitoring} without notification. For the case the TTC drops further and equation ~\eqref{eq:TimeToEvade} is satisfied, then the state switches to \textit{AES in regulation}. In this state the evasive manoeuvre is initiated. If the current path becomes infeasible, a new path can be selected by the replanning function. When the manoeuvre is fully completed, AES will go back to the AES monitoring state. In the case that a failure is detected (e.g., the replanning function cannot provide a feasible path or a system malfunction prevents successful completion) the system will go into an \textit{AES aborted} state. Depending on the reason for abortion, the system can be reinitialized to AES Standby or will be turned off. 

\subsection{Motion control}

The purpose of the lateral controller is to make the vehicle follow a predefined evasive path. Where normally a steering action is handled by the steering system, a combination of steering and individual wheel braking can enhance the vehicle steering capabilities. 
For an evasive manoeuvre, it is beneficial to be able to act as late as possible. The brake system can be helpful here as it can produce a yaw moment more quickly. Furthermore, in a driving assistance system, the driver has to be taken into account. Unless the vehicle is driving autonomously at SAE level 3 or higher, the driver is expected to keep his hands on the wheel. functional safety studies (out of the scope of this article) have shown that the best steering performances under the constraint of driver controllability in case of failure are achieved with a combination of steering and differential braking, rather than steering only. Individual wheel braking can help here to compensate for the driver limitation.

\begin{figure}[t]
\fontsize{7pt}{7pt}\selectfont
\def\svgwidth{\columnwidth}
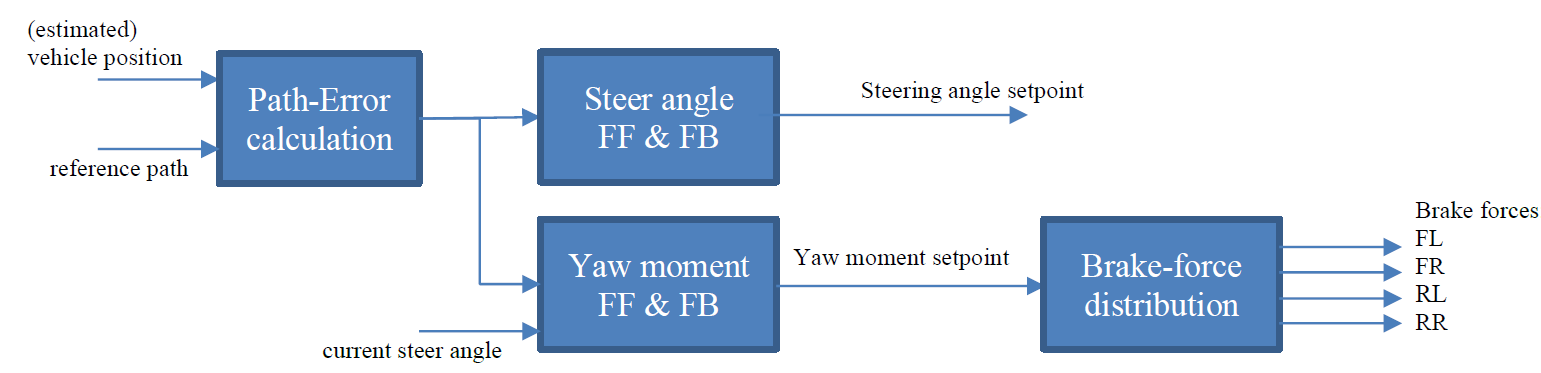
    \caption{Sketch of AES control architecture.}
    \label{fig:AESControlStructure}
\end{figure}

In Fig.~\ref{fig:AESControlStructure} the controller architecture is shown. The control part is split up into a steering controller and a yaw moment controller. The inputs of the steering controller are based on the path tracking errors. The yaw moment controller uses the path tracking errors and the actual steering angle as an input to generate the required brake forces. By using the actual steering angle, the driver input is incorporated as a disturbance, enabling the controller to work around the driver.

The controller brings the vehicle to the intended path by regulating the path-tracking errors to zero. The path-tracking errors are described by the lateral error $y_e$ and the heading error $\psi_e$.
 To calculate these path-tracking errors, first the path needs to be transformed from the local frame (the vehicle frame at the AES trigger) to the current vehicle frame. This is done by multiplying the path by the inverse transformation matrix as follows:
\begin{equation}
\left[\begin{matrix}x_{path}\\y_{path}\\1\\\end{matrix}\right]=\left[\begin{matrix}\cos{\left(\psi\right)}&-\sin{\left(\psi\right)}&X\\\sin{\left(\psi\right)}&\cos{\left(\psi\right)}&Y\\0&0&1\\\end{matrix}\right]^{-1}\left[\begin{matrix}X_{path}\\Y_{path}\\1\\\end{matrix}\right]
\end{equation}

\begin{figure}[t]
\includegraphics[width=\columnwidth]{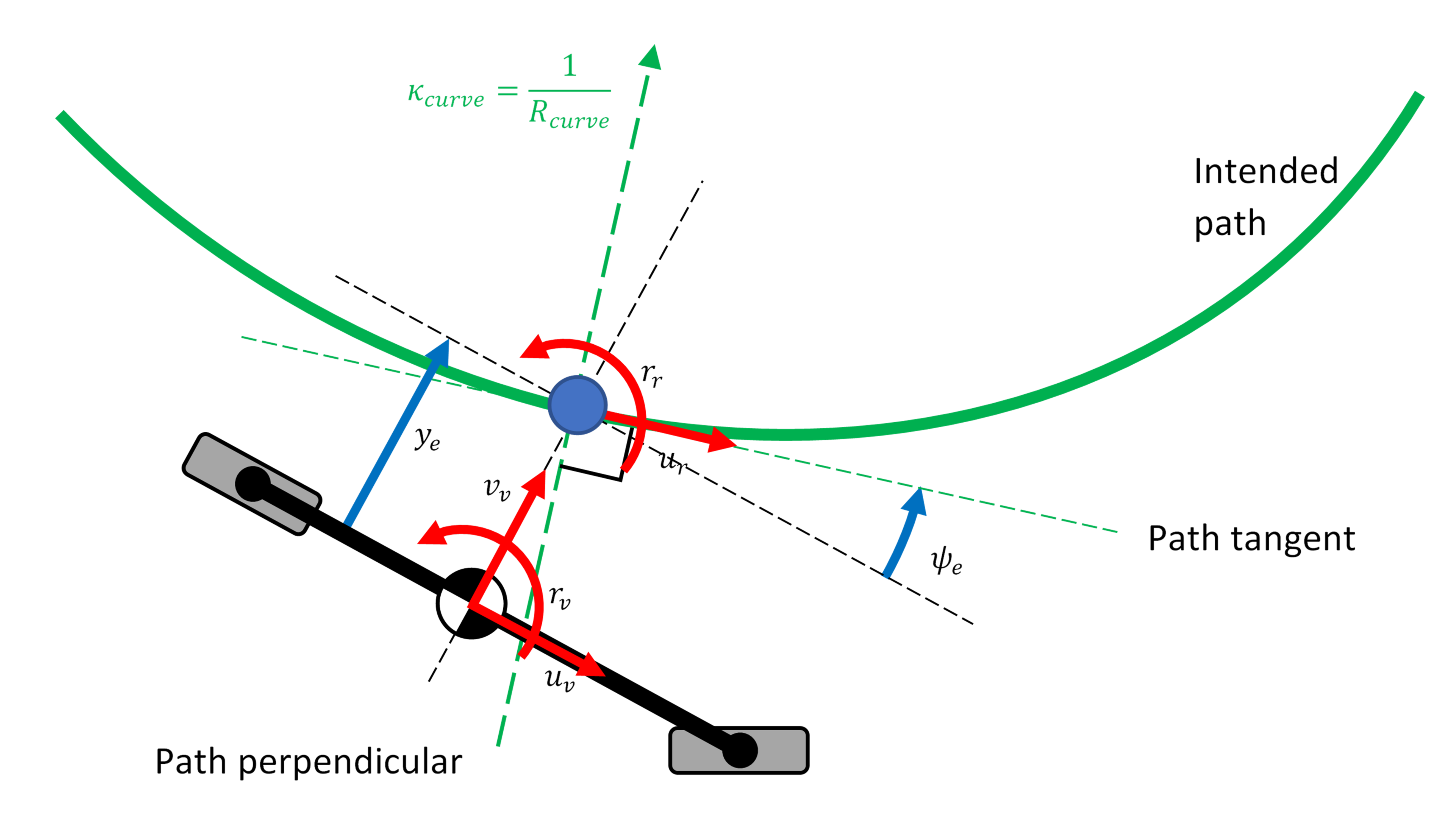}
\caption{Path error calculation}
\label{fig:PathErrorCal}
\end{figure}

\begin{figure}[t]
\includegraphics[width=\columnwidth]{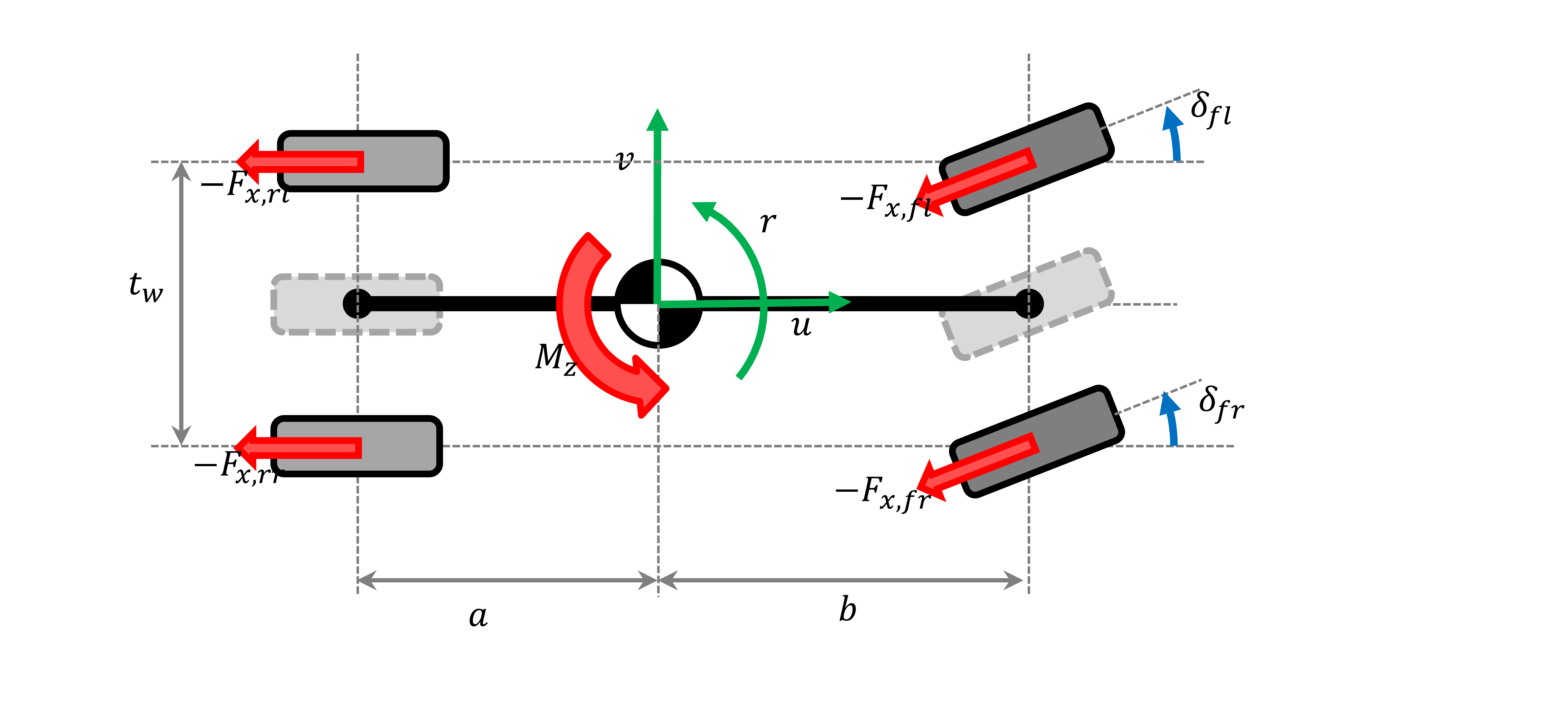}
\caption{Concept of differential braking}
\label{fig:ConceptDiffBraking}
\end{figure}

In this equation $X$,$Y$ and $\psi$ are the location and heading of the ego vehicle in the local frame, $X_{path}$,$Y_{path}$ are the path coordinates in the local frame and $x_{path}$ and $y_{path}$ are the path coordinates in the vehicle’s frame. After transformation the path-tracking errors can then be calculated as depicted in Fig.~\ref{fig:PathErrorCal}. \\The controlled vehicle is modeled using the commonly known linear single-track model, augmented with an additional external yaw-moment input for the purpose of differential braking. The model is described by the following equations of motion (as also depicted in Fig.~\ref{fig:PathErrorCal})
\begin{align}
    \begin{bmatrix}
        \dot{v}_v\\\dot{r}    
    \end{bmatrix}=&
    \begin{bmatrix}
        \frac{C_f+C_r}{mu_v}&\frac{aC_f-bC_r}{mu_v}-u_v\\\frac{aC_f-bC_r}{I_{zz} u_v}&\frac{a^2C_f+b^2C_r}{I_{zz}u_v}
    \end{bmatrix}
    \begin{bmatrix}
        v_v\\r  
    \end{bmatrix}\nonumber\\&+
    \begin{bmatrix}
        -\frac{C_f}{m}&0\\
        -\frac{aC_f}{I_{zz}}&\frac{1}{I_{zz}}
    \end{bmatrix}
    \begin{bmatrix}
        \delta_g\\M_{z,ext}
    \end{bmatrix},\label{eq:statespace}
\end{align}
where, $I_{zz}$ represents the moment of inertia of the vehicle around the vertical axis, $r$ represents the yaw-rate, $v_v$ represents the lateral velocity, $u_v$ represents the longitudinal velocity and $M_{z,ext}$ represents the externally applied yaw-moment, which is longitudinally induced by the tyre braking forces as follows
\begin{align}
    M_{z,ext}=\frac{1}{2}t_w(-F_{x,fl}+F_{x,fr}-F_{x,rl}+F_{x,rr})    
\end{align}
where the trackwidth $t_w$ and the tyre forces $F_{x,fl},F_{x,fr},F_{x,rl},F_{x,rr}$ are depicted in Fig.~\ref{fig:ConceptDiffBraking}. The vehicle dynamic poles, i.e., the eigenvalues of the state matrix of~\eqref{eq:statespace}, are from now on denoted as $\sigma_1^v,\sigma_2^v$. 
 
With some mathematical manipulation, the following model is used to describe the error’s dynamics of the vehicle position w.r.t the path to follow:
\begin{equation}
    v_v=u_v\psi-{\dot{y}}_e\\
    r_v=u_v\kappa-{\dot{\psi}}_e
\end{equation}

Finally, the augmented bicycle model with both steering and differential braking is presented as follows:
\begin{align}
    \footnotesize
    \setlength{\arraycolsep}{2.5pt}
    \medmuskip = 1mu 
    &\begin{bmatrix}{\dot{y}}_e\\{\ddot{y}}_e\\{\dot{\psi}}_e\\{\ddot{\psi}}_e\\\end{bmatrix}=\begin{bmatrix}0&1&0&0\\0&\frac{C_{f}+C_{r}}{mu}&-\frac{C_{f}+C_{r}}{m}&\frac{aC_{f}-bC_{r}}{mu}\\0&0&0&1\\0&\frac{aC_{f}-bC_{r}}{I_{zz}u}&-\frac{aC_{f}-bC_{r}}{I_{zz}}&\frac{a^2C_{f}+b^2C_{r}}{I_{zz}u}\\\end{bmatrix}\begin{bmatrix}y_e\\{\dot{y}}_e\\\psi_e\\{\dot{\psi}}_e\\\end{bmatrix}\nonumber\\&+\begin{bmatrix}0&0\\\frac{C_{f}}{m}&0\\0&0\\\frac{aC_{f}}{I_{zz}}&\frac{-1}{I_{zz}}\\\end{bmatrix}\!\begin{bmatrix}\delta_g\\M_{z,ext}\\\end{bmatrix}\!+\!\begin{bmatrix}0&0\\u_v^2-\frac{aC_{f}-bC_{r}}{m}&0\\0&0\\-\frac{a^2C_{f}+b^2C_{r}}{I_{zz}}&u_v\\\end{bmatrix}\!\begin{bmatrix}\kappa\\\dot{\kappa}\\\end{bmatrix}
\end{align}

\subsubsection{Mixed lateral controller}
The proposed controller composes of a state-feedback and a feedforward term. This choice is linked to the embeddability requirement in ECUs and the tunning’s simplicity
\begin{equation}
    u=Ky+u_{ff}
\end{equation}
where
\begin{equation}
u_{ff}=\left(\begin{matrix}\delta_g^{ff}\\M_{z,ext}^{ff}\\\end{matrix}\right).
\end{equation}

The controller of the matrix is defined as
\begin{equation}
    K=\left[\begin{matrix}k_y^\delta&k_{Dy}^\delta&k_\psi^\delta&k_{D\psi}^\delta\\k_y^M&k_{Dy}^M&k_\psi^M&k_{D\psi}^M\\\end{matrix}\right]
\end{equation}

The feedforward steer angle term ($\delta_g^{ff}$) and the feedforward yaw moment ($M_{z,ext}^{ff}$) can be found by looking at the steady-state (circle) solution of the vehicle dynamics model, with $r_v=u_v \kappa$  (assuming small body slip angles). 
For steering control only ($M_{z,ext}^{ff}=0$) the feedforward steering angle becomes: 
\begin{equation}
\delta_g^{ff}=\left(K_\delta u_v^2+l_{wb}\right)\kappa
\end{equation}
Where $K_\delta$ is the understeer coefficient, defined as: 

\begin{equation}
    K_\delta=\frac{m}{l_{wb}}\left(\frac{a}{C_{r}}-\frac{b}{C_{f}}\right)
\end{equation}

For differential braking control only, the steering angle can be seen as a driver disturbance. This means that the steering angle can take on any value. The feedforward yaw moment therefore becomes: 

\begin{equation}
    M_{z,ext}^{ff}=\frac{l_{wb}C_{f}C_{r}}{C_{f}{+C}_{\alpha r}}\left(\delta_g-\left(K_\delta u_v^2+l_{wb}\right)\kappa\right)
\end{equation}

Now that the feedforward terms are known, the feedback solution can be determined. 
Solving the characteristic equation leads to the four poles of the system ($\sigma_i$). As already mentioned in the derivation of the vehicle model, the vehicle model itself has 2 poles ($\sigma_{1}^v$ and $\sigma_{2}^v$). Our aim is to not change the location of these two vehicle poles. The remaining two poles of the closed loop system ($\sigma_1$ and $\sigma_2$) can however be chosen freely. The resulting characteristic equation for this desired closed loop system now becomes:
\begin{equation}
\label{eq:chareq}
\left(\lambda-\sigma_1\right)\left(\lambda-\sigma_2\right)\left(\lambda-\sigma_1^v\right)\left(\lambda-\sigma_2^v\right)=0    
\end{equation}

\subsubsection{Feedback steering control} 
The characteristic equation~\eqref{eq:chareq} is solved for the steering angle control gains, through preservation of the vehicle poles, $\sigma_1^v$ and $\sigma_2^v$, which leads to the solution:

\begin{equation}
    k_y^\delta=\left(\frac{l_{wb}}{u^2}+\frac{m\left(aC_{f}-bC_{r}\right)}{l_{wb}C_{f}C_{r}}\right)\sigma_1\sigma_2
\end{equation}

\begin{equation}
    k_{Dy}^\delta=-\frac{m\left(aC_{f}-bC_{r}\right)}{l_{wb}C_{f}C_{r}}\left(\sigma_1+\sigma_2\right)-\frac{ma}{C_{r}u}\sigma_1\sigma_2
\end{equation}

\begin{equation}
    k_\psi^\delta=-\frac{l_{wb}}{u}\left(\sigma_1+\sigma_2\right)-\left(\frac{\left(C_{f}+C_{r}\right)I_{zz}}{C_{f}C_{r}l_{wb}}+\frac{{bl}_{wb}}{u^2}\right)\sigma_1\sigma_2
\end{equation}

\begin{equation}
    k_{D\psi}^\delta=\frac{\left(C_{f}+C_{r}\right)I_{zz}}{C_{f}C_{r}l_{wb}}\left(\sigma_1+\sigma_2\right)+\frac{I_{zz}}{C_{r}u}\sigma_1\sigma_2
\end{equation}

\subsubsection{Feedback differential braking control}
Similar to the steering controller, the equation~\eqref{eq:chareq} is solved for differential braking control gains, which leads to the solution:

\begin{align}
    k_y^M=&\left(\frac{C_{f}C_{r}l_{wb}^2}{u^2\left(C_{f}+C_{r}\right)}+\frac{m\left(aC_{f}-bC_{r}\right)}{\left(C_{f}+C_{r}\right)}\right)\sigma_1\sigma_2\\
    k_{Dy}^M=&-\frac{mC_{f}C_{r}l_{wb}^2}{{u\left(C_{f}+C_{r}\right)}^2}\sigma_1\sigma_2\nonumber\\&+\left(\frac{ml_{wb}C_{f}}{C_{f}+C_{r}}-mb\right)\left(\sigma_1+\sigma_2\right)
\end{align}

\begin{align}
    k_\psi^M=&\frac{C_{f}C_{r}l_{wb}^2}{u\left(C_{f}+C_{r}\right)}\left(\sigma_1+\sigma_2\right)\nonumber\\&+\left(I_{zz}-\frac{C_{f}C_{r}l_{wb}^2\left(aC_{f}-bC_{r}\right)}{{u^2\left(C_{f}+C_{r}\right)}^2}\right)\sigma_1\sigma_2
\end{align}
\begin{equation}
    k_{D\psi}^M=-I_{zz}\left(\sigma_1+\sigma_2\right)
\end{equation}
Finally, an uncovered topic is the translation of a certain closed-loop control signal $M_{z.ext}$ to the separate brake forces we actually have access to. This is described through the following allocation equations
\begin{align}
    F_{brake,fl}=&\max(0,\frac{t_wM_{z,ext}}{2})\cdot i_f\\
    F_{brake,fr}=&-\min(0,\frac{t_wM_{z,ext}}{2})\cdot i_f\\
    F_{brake,rl}=&\max(0,\frac{t_wM_{z,ext}}{2})\cdot i_r\\
    F_{brake,rr}=&-\min(0,\frac{t_wM_{z,ext}}{2})\cdot i_r
\end{align}
where $i_{f}$ is the factor of brake-force applied to the front wheels and $i_{r}$ is the factor of brake-force applied to the rear wheels. These factors are tuning factors which depend on the particular vehicle dynamical characteristics and are chosen such that $i_{f}+i_{r} = 1$.
\section{Simulation and Experiments}
\begin{figure}[t]
    \centering
    \includegraphics[width=\textwidth]{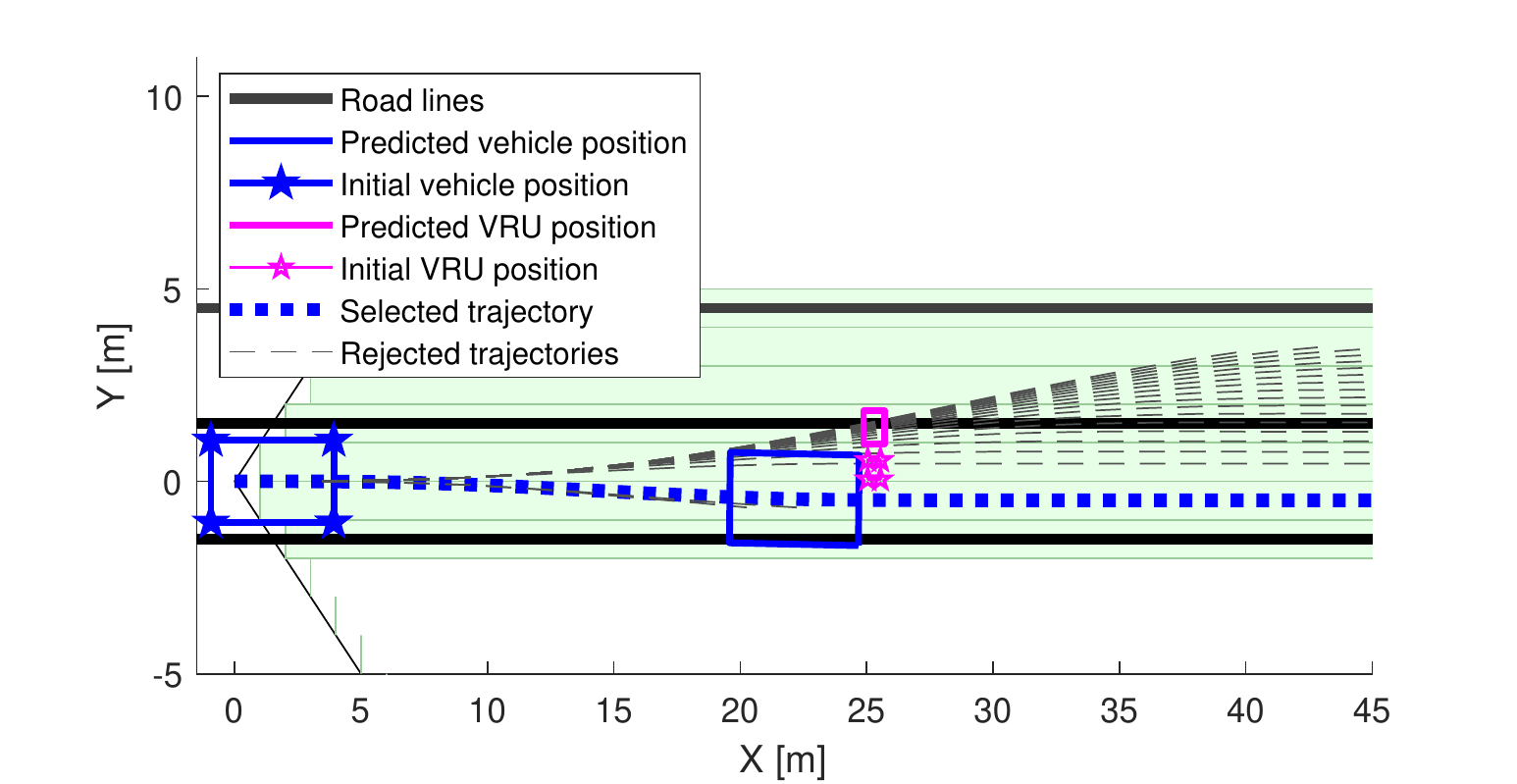}
    \caption{Proposed trajectories, including the selected and rejected trajectories, for avoiding a collision with the VRU.}
    \label{fig:sim1}
\end{figure}

\begin{figure}[t]
    \centering
    \includegraphics[width=\textwidth]{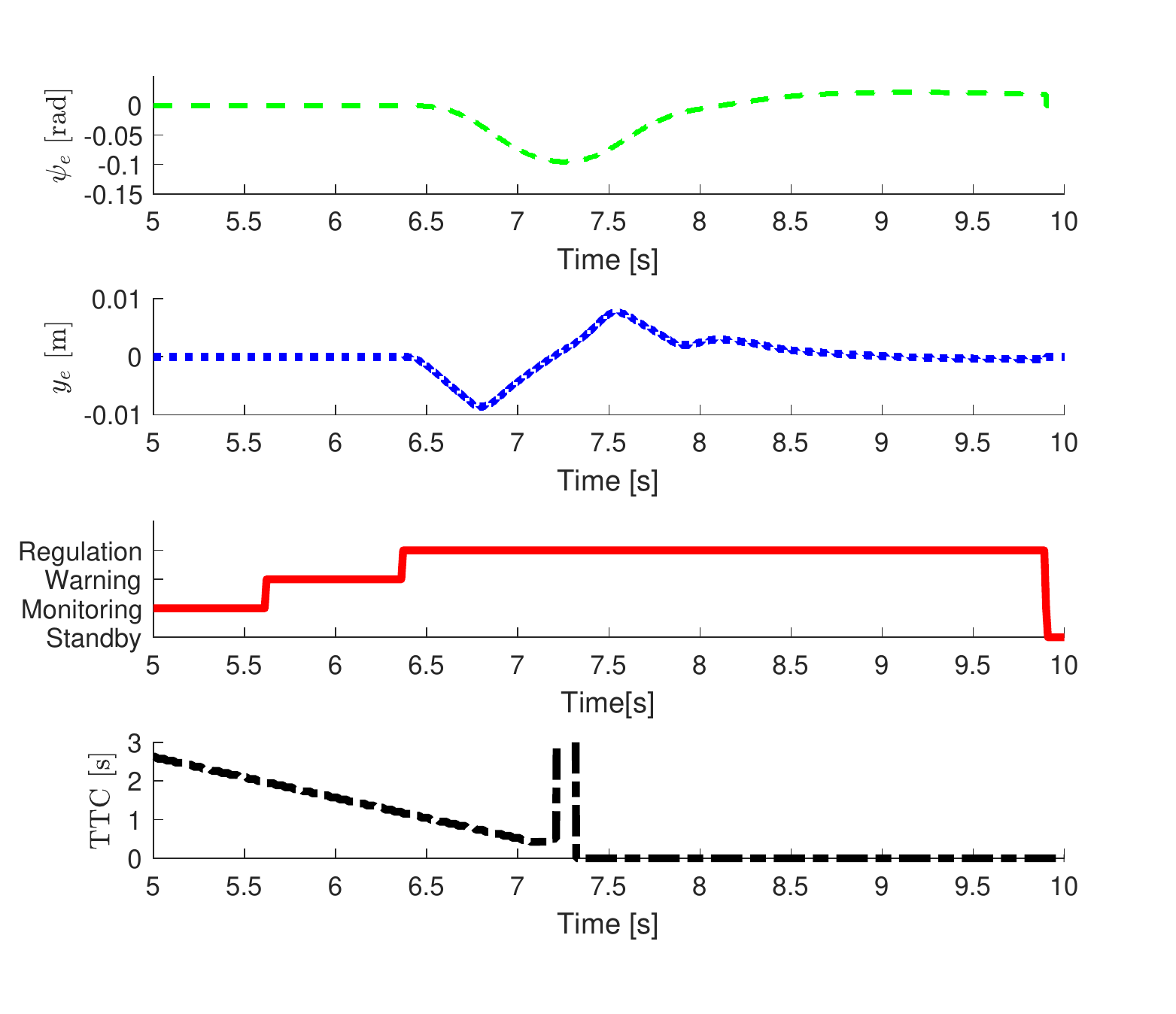}
    \caption{States $y_e,\:\psi_e$ relative to the AES path, the AES supervisory state and the TTC relative to the VRU. }
    \label{fig:sim2}
\end{figure}

\begin{figure}[t]
    \centering
    \includegraphics[width=\textwidth]{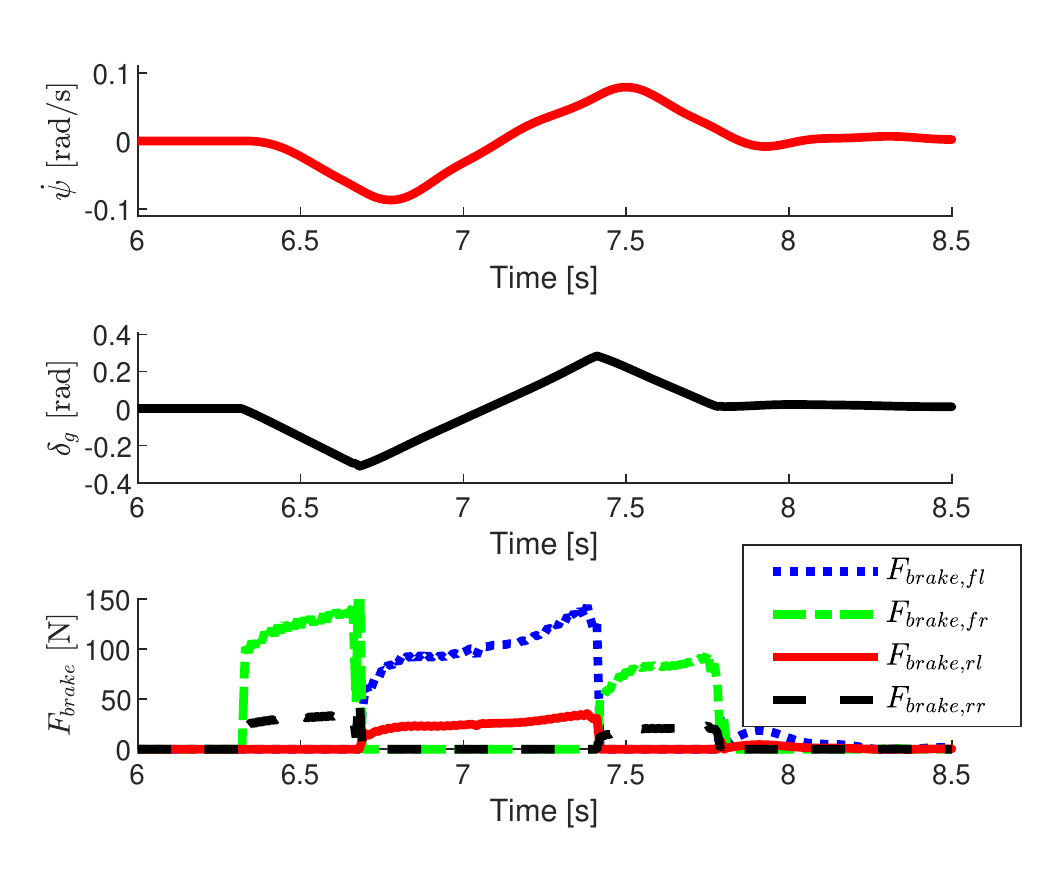}
    \caption{State $\dot{\psi}$, the ground steering angle $\delta_g$ and the braking torques for all four individual wheels.}
    \label{fig:sim3}
\end{figure}

In this section, the effectiveness and applicability of the proposed approach is demonstrated through a simulation study and real-life experiments. In the simulation study, we consider the use-case which has been proposed in section~3, Figure~\ref{fig:UCselection}, where a VRU crosses the road at a velocity of $1\:\text{m}\cdot\text{s}^{-1}$. The model, used to simulate the ego vehicle, is a double-track model with independently controlled brake actuators, identified from the experimental vehicle. A Magic Formula (Pacejka) tyre model is employed to characterize the tyre dynamics. This high level of fidelity is used to capture the behavior of the vehicle in highly dynamic situations like an AES manoeuvre. The scenario is depicted in Figure~\ref{fig:sim1}. The ego vehicle drives at a longitudinal velocity of $20\:\text{m}\cdot \text{s}^{-1}$ The boundaries of the road form the extrema for the driveable space, as such, the generated trajectories are not allowed to surpass this enclosed space. Figure~\ref{fig:sim1} gives a planar representation of the result. It can be observed, that the most effective and last possible AES manuever, at a TTC of $0.42s$, is described by passing the crossing pedestrian on the right. From $7.2s$ onwards, the AES system has avoided the object and, hence, the TTC is tends
towards infinity, since the object is no longer on the longitudinal path of the ego vehicle. Note, that for conventional AEB systems, assuming a deceleration of $11\:\text{m}\cdot \text{s}^{-2}$ it would have been impossible to avoid such a collision. The state of the vehicle, with respect to the path, the VRU and of the AES supervisory control, is depicted in Figure~\ref{fig:sim2}. The yaw-rate of the vehicle and the actual applied control inputs are depicted in Fig.~\ref{fig:sim3}. It can be observed that the lateral deviation with respect to the path lies within a margin of $0.01$m.

\begin{figure}[t]
    \centering
    \includegraphics[width=\columnwidth]{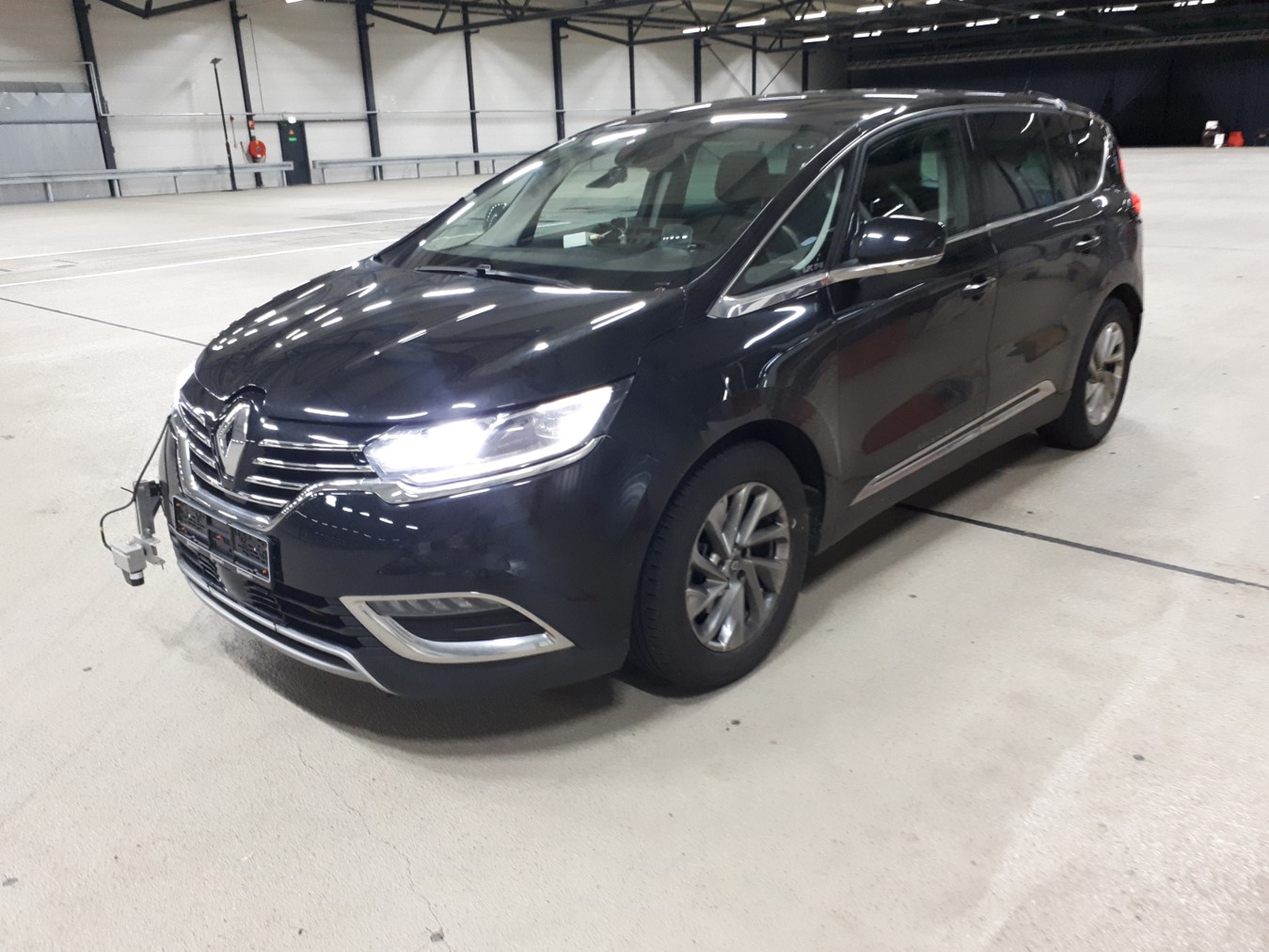}
    \caption{Testing platform for the experimental results.}
    \label{fig:espace}
\end{figure}

To globally validate the developed control approach, some experiments are conducted on a Renault Espace testing vehicle~\ref{fig:espace}. The vehicle is equipped with a limited sensor set and is able to receive a steering wheel setpoint and an independent brake torque setpoint on each wheel. For safety reasons, the actuation capabilities are limited, i.e., the driver always needs to be able to overrule or take over control from the system. This will also be the case in a commercial application. The exact values and how these limits are derived will not be elaborated. Moreover it needs to be mentioned that the motion controller was not extensively tuned since it concerns a preliminary feasibility study. The main goal of these experiments is to validate that the evasive path planned, can be realized on a global level. Experiments are conducted at a driving speed of 65 km/h using steering only, steering in combination with braking and braking only to perform an evasive manoeuvre. The reference path consists of a single lane change to avoid a collision with an object that appears in front of the vehicle. 

\begin{figure}[t]
\includegraphics[width=\columnwidth]{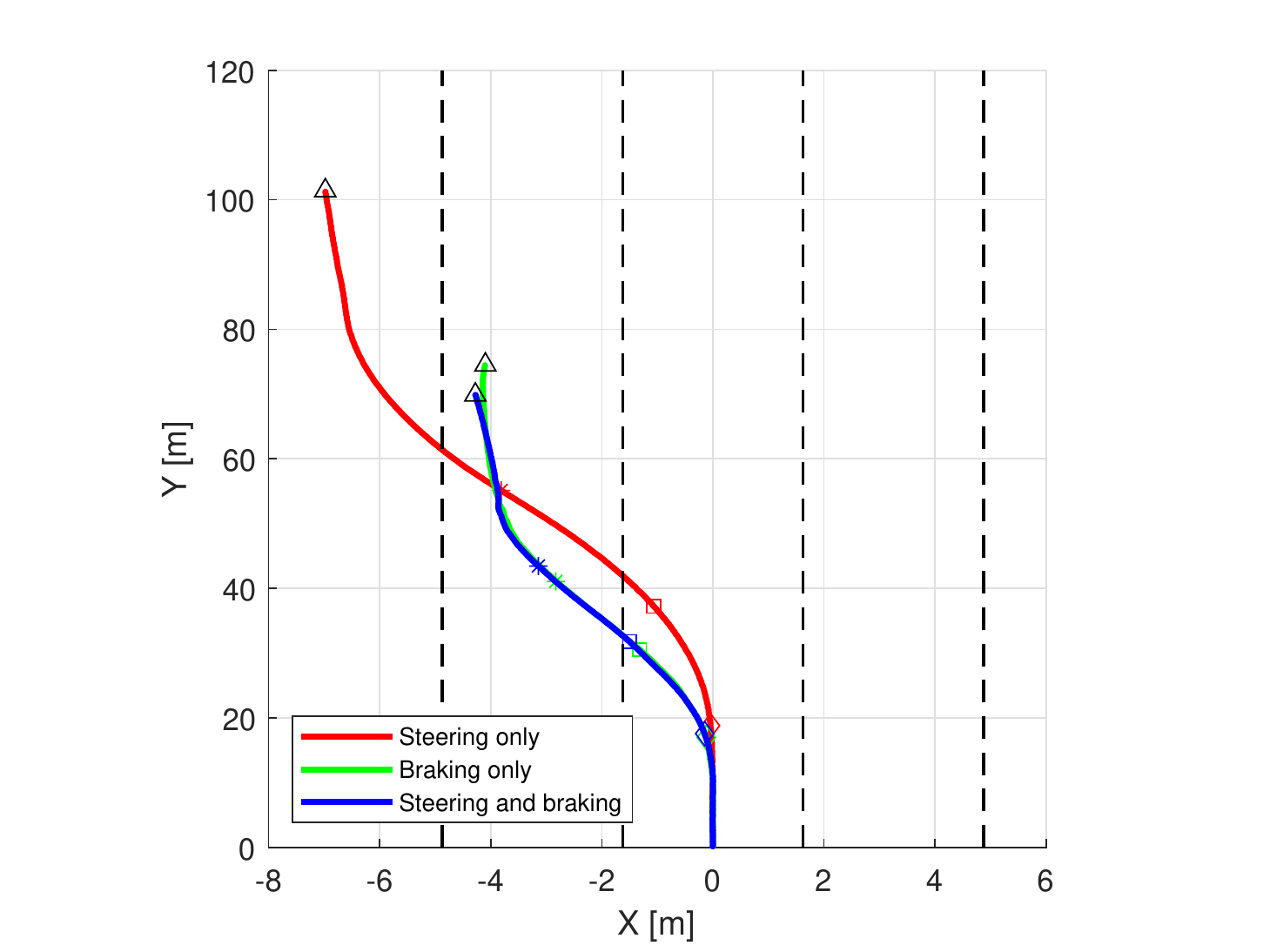}
\caption{Trace of the vehicle's rear axle center point for an evasive manoeuvre performed with three control approaches: steering only, braking only, steering and braking combined.}
\label{fig:ExperimentsFigure1}
\end{figure}

In Fig. \ref{fig:ExperimentsFigure1}, the results of three experiments are shown in an X-Y graph. On the vertical axis and horizontal axis respectively the longitudinal and lateral travelled distance is plotted. The four vertical dashed lines indicate the lane edges of a typical European lane width of 3.25 meter. In the origin of the plot the actuation during three experiments is initiated; steering only, braking only and steering and braking. The GPS traces of the vehicle's rear axle center are plotted with a dashed, dash-dot and continuous line respectively. 

In order to gain insight in the lateral displacement over time a diamond, square and star shaped marker is added respectively at one, two and three seconds after the initiation of the manoeuvre. The triangular markers indicate the end of the manoeuvre. As can be seen, the steering-only approach builds up the lateral displacement significantly later as compared to the brake-only approach. Moreover poor controller performance results in a large overshoot and thus a shift of two lanes to the left. As can be detected from the placement of the diamond, square and star markers of braking and braking combined with steering, the braking combined with steering builds the same lateral displacement up slightly faster. The difference in lateral displacement build up between steering only and braking only is mainly caused by the phase shift in the response of the steering system. 
The results show, that when the brakes are used, the fastest response can be realized. 

\begin{figure}[t]
\includegraphics[width=1.05\columnwidth]{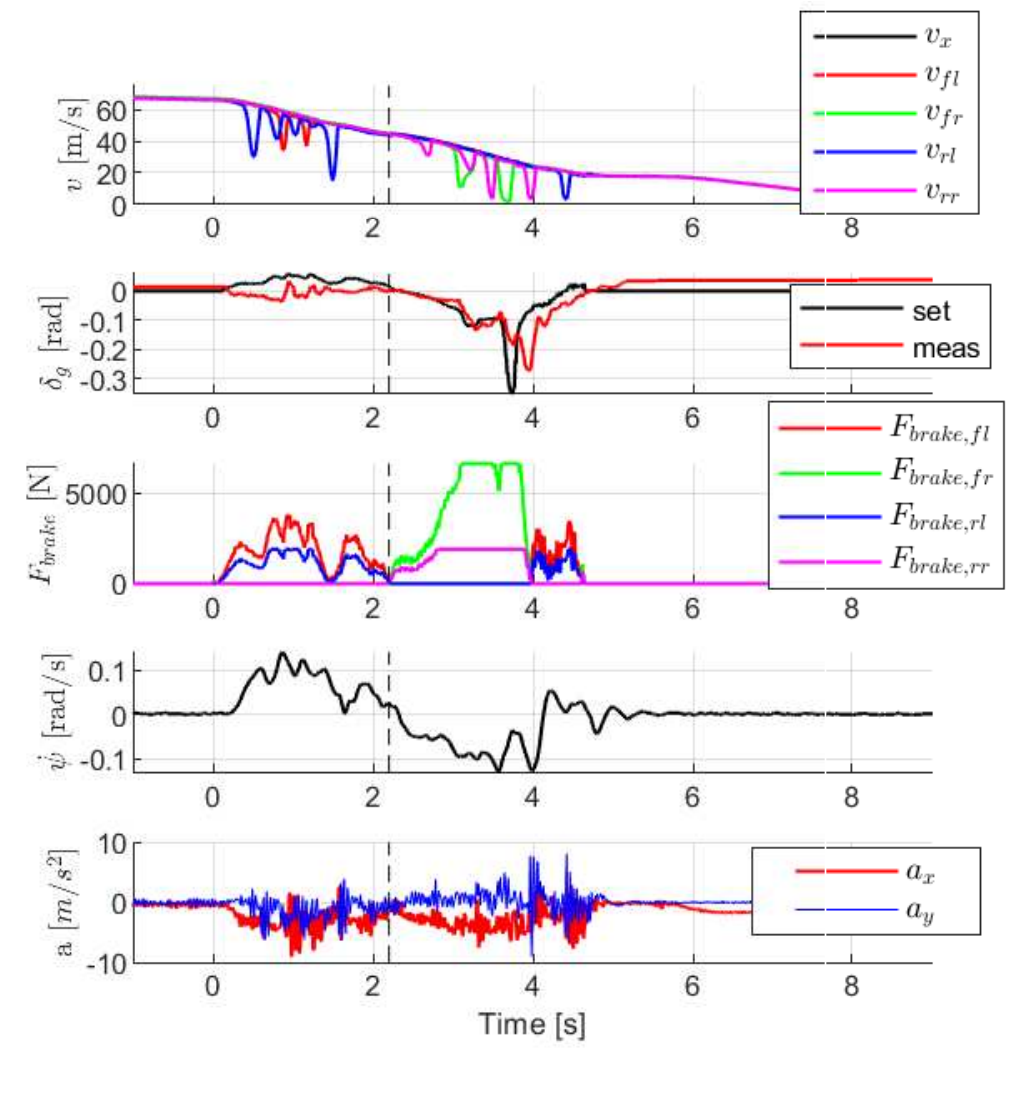}
\caption{Measurement results of the evasive manoeuvre performed using the combination of steering and braking.}
\label{fig:ExperimentsFigure2}
\end{figure}

In Fig. \ref{fig:ExperimentsFigure2} different signals are depicted of the test in which the combination of braking and steering is used to evade. In the top plot the vehicle and independent wheel speeds are shown. It shows the decrease in speed generated by the use of the brakes during the evasive manoeuvre that is initiated at t=0sec. It also shows that individual wheels tend to lock which is caused by the limited performance of the prototype ABS. From $t=5.8s$ on, the driver applies the brakes to bring the vehicle to standstill. 

In the second plot, the steering wheel control setpoint and measured value are shown. What stands out is that despite a slight positive setpoint the measured value is even in the opposite negative direction. This behavior is caused by the vehicle’s steering system and suspension (caster and toe)  geometry. For the vehicle used it means that asymmetric braking at the front axle will result in the steering wheel rotating towards the wheel that is braked the least e.g. the steering system actually counteracts the lateral displacement generated with the brakes. Moreover it can be seen that the controller uses more steering action to stabilize the vehicle in the adjacent lane at the end of the maneouvre. From actuation point of view, counter steering is initiated at $2.2s$, as indicated by the dashed vertical line in the subplots. 

In the third plot, the setpoints for the individual brakes are shown. We can also see here clearly that, till 2.2 s, the front-left and rear-left brakes are used while, after that time, mainly the front-right and rear-right brakes are used to stabilize the vehicle. In the subplot also the saturation level can be detected from which it can be concluded that the controller is not fully utilizing the actuation possibilities between $0s$ and $0.8 s$. In the two lower plots, the yawrate and lateral and longitudinal accelerations are shown. These show that the vehicle decelerates with approximately $4\: \text{m}\cdot \text{s}^{-2}$ during the manoeuvre. 

From the tests, it can be concluded that an evasive manoeuvre can be performed by using steering, differential braking or a hybrid combination. Differential braking mitigates the phase lag of the steering system and can reduce the steering disturbance induced by the driver (for example, the driver tightly holding the steering wheel). As a matter of fact, it appears that differential braking could be used as a standalone solution for AES. However, it is clear that the total lateral capability is increased when using the combination; moreover, using the steering wheel to actuate improves the driver involvement, which is essential from a driver-acceptance standpoint. Finally, the experimental results show that the software pipeline runs on a real vehicle.

\section{Conclusion}
In this work we have provided a complete architecture for the design of an AES system as a safety feature for vehicles. After introducing the state of the art, the technical details of our architecture have been provided, including a novel method for online path-replanning, as well as a novel methodology for defining an AES trigger. The proposed approach has been proven both in simulation and partially in experiments, showing that the architecture functions, both in simulation and in real-life. 
    \printbibliography

@article{Jonasson2017,
   author = {M. Jonasson and M. Thor},
   journal = {Vehicle System Dynamics},
   month = {5},
   pages = {791-809},
   publisher = {Taylor & Francis},
   title = {Steering redundancy for self-driving vehicles using differential braking},
   volume = {56},
   year = {2017},
}

@article{Breuer1998,
    author = {J. Breuer},
    journal = {Int. Technical Conference on ESV},
    title = {Analysis of driver-vehicle-interactions in an evasive manoeuvre-results of moose test studies},
    year = {1998},
}

@article{Hayashi:2012,
author = { Ryuzo   Hayashi  and  Juzo   Isogai  and  Pongsathorn   Raksincharoensak  and  Masao   Nagai },
title = {Autonomous collision avoidance system by combined control of steering and braking using geometrically optimised vehicular trajectory},
journal = {Vehicle System Dynamics},
volume = {50},
number = {sup1},
pages = {151-168},
year  = {2012},
publisher = {Taylor & Francis},
}

@ARTICLE{Brannstrom:14,
author={Mattias {Brännström} and Erik {Coelingh} and  J {Sjoberg}},
journal={Int. J. Vehicle Safety},
title={Decision-making on when to brake and when to steer
to avoid a collision},
year={2014},
month = {01},
pages = {87-106},
volume={7},
number={1},}

@article{Ackermann:14,
title = {Collision avoidance with automatic braking and swerving},
journal = {19th IFAC World Congress},
volume = {47},
number = {3},
pages = {10694-10699},
year = {2014},
author = {Carlo Ackermann and Rolf Isermann and Sukki Min and Changwon Kim},
}

@misc{Sim:17,
 title     = "LANE CHANGE CONTROL DEVICE AND CONTROL METHOD ({P}atent)",
 number    = "20160185388",
 author    = "Sim, Sang Kyun ",
 year      = "2016",
}

@misc{Arbitmann:13,
	author={Maxim {Arbitmann} and Matthias {Schorn} and
Rolf {Isermann}},
	title={Method and device for steering a motor vehicle ({P}atent)},
	language={English},
	assignee={ Continental Teves AG & Co. oHG},
	address={Pyeongtaek-si, Gyeonggi-do (KR)},
	nationality ={German},
	type={US Patent},
	note={8355842 B2},
	year={2013},
}

@misc{Hong:2017,
	author={Sun {Hong} and Zhang {Changzhu} and Chen {Qijun} and Shen {Mengjiao} and Chen {Longquan} and An {Guangyong}},
	title={State feedback based vehicle path tracking {H}$\infty$ control method ({P}atent)},
	language={Chinese},
	assignee={Tongji University},
	nationality ={Chinese},
	type={China Patent},
	note={107015477A},
	year={2017},
}

@ARTICLE{Funke:17,
author={J. {Funke} and M. {Brown} and S. M. {Erlien} and J. C. {Gerdes}},
journal={IEEE Transactions on Control Systems Technology},
title={Collision Avoidance and Stabilization for Autonomous Vehicles in Emergency Scenarios},
year={2017},
volume={25},
number={4},
pages={1204-1216},}

@article{Berntorp:17,
author={K. {Berntorp}},
journal={2017 American Control Conference (ACC)},
title={Path planning and integrated collision avoidance for autonomous vehicles},
year={2017},
volume={},
number={},
pages={4023-4028},}

@inproceedings{serafimguardini:22,
  TITLE = {{Employing Severity of Injury to Contextualize Complex Risk Mitigation Scenarios}},
  AUTHOR = {Serafim Guardini, Luiz Alberto and Spalanzani, Anne and Laugier, Christian and Martinet, Philippe and Do, Anh-Lam and Hermitte, Thierry},
  BOOKTITLE = {{IV 2020 -- 31st IEEE Intelligent Vehicles Symposium}},
  ADDRESS = {Las Vegas / Virtual, United States},
  PUBLISHER = {{IEEE}},
  PAGES = {1-7},
  YEAR = {2020},

}

@inproceedings{Giugliano2015,
  year = {2015},
  month = {10},
  publisher = {American Society of Mechanical Engineers},
  author = {Luke Giugliano and Craig E. Beal},
  title = {Dynamic Rear-End Collision Mitigation for a Vehicle About to be Struck},
  booktitle = {ASME 2015 Dynamic Systems and Control Conference}
}

@article{Gallen2013,
  title={Supporting Drivers in Keeping Safe Speed in Adverse Weather Conditions by Mitigating the Risk Level},
  author={Romain Gallen and Nicolas Hauti{\`e}re and Aur{\'e}lien Cord and Sebastien Glaser},
  journal={IEEE Transactions on Intelligent Transportation Systems},
  year={2013},
  volume={14},
  pages={1558-1571}, }

@article{Lee2019,
  title={Collision Avoidance/Mitigation System: Motion Planning of Autonomous Vehicle via Predictive Occupancy Map},
  author={Kibeom Lee and Dongsuk Kum},
  journal={IEEE Access},
  year={2019},
  volume={7},
  pages={52846-52857}, }

@INPROCEEDINGS{Lima:2015,
  author={Lima, Pedro F. and Trincavelli, Marco and Martensson, Jonas and Wahlberg, Bo},
  booktitle={2015 European Control Conference (ECC)}, 
  title={Clothoid-based model predictive control for autonomous driving}, 
  year={2015},
  volume={},
  number={},
  pages={2983-2990},}

@INPROCEEDINGS{Williams:2016,  author={Williams, Grady and Drews, Paul and Goldfain, Brian and Rehg, James M. and Theodorou, Evangelos A.},  booktitle={2016 IEEE International Conference on Robotics and Automation (ICRA)},   title={Aggressive driving with model predictive path integral control},   year={2016},  volume={},  number={},  pages={1433-1440}}

@article{DO:2021,
title = {LPV approach for collision avoidance: Controller design and experiments},
journal = {Control Engineering Practice},
volume = {113},
pages = {104856},
year = {2021},
author = {Anh-Lam Do and François Fauvel},
}

@article{Keller:2011,
author = {Keller, Christoph and Dang, Thao and Fritz, Hans and Joos, Armin and Rabe, Clemens and Gavrila, Dariu},
year = {2011},
month = {05},
pages = {1292 - 1304},
title = {Active Pedestrian Safety by Automatic Braking and Evasive Steering},
volume = {12},
journal = {Intelligent Transportation Systems, IEEE Transactions on}
}

@incollection{Dang2012,
author = {Dang, Thao and Desens, Jens and Franke, Uwe and Gavrila, Dariu and Schäfers, Lorenz and Ziegler, Walter},
year = {2012},
month = {01},
pages = {759-782},
title = {Steering and Evasion Assist},
booktitle = {Handbook of Intelligent Vehicles},
publisher="Springer London",
isbn = {978-0-85729-084-7},
}

@ARTICLE{Gonzalez:2016,
  author={González, David and Pérez, Joshué and Milanés, Vicente and Nashashibi, Fawzi},
  journal={IEEE Transactions on Intelligent Transportation Systems}, 
  title={A Review of Motion Planning Techniques for Automated Vehicles}, 
  year={2016},
  volume={17},
  number={4},
  pages={1135-1145},}

@ARTICLE{Ziegler:2014,
  author={Ziegler, Julius and Bender, Philipp and Schreiber, Markus and Lategahn, Henning and Strauss, Tobias and Stiller, Christoph and Dang, Thao and Franke, Uwe and Appenrodt, Nils and Keller, Christoph G. and Kaus, Eberhard and Herrtwich, Ralf G. and Rabe, Clemens and Pfeiffer, David and Lindner, Frank and Stein, Fridtjof and Erbs, Friedrich and Enzweiler, Markus and Knöppel, Carsten and Hipp, Jochen and Haueis, Martin and Trepte, Maximilian and Brenk, Carsten and Tamke, Andreas and Ghanaat, Mohammad and Braun, Markus and Joos, Armin and Fritz, Hans and Mock, Horst and Hein, Martin and Zeeb, Eberhard},
  journal={IEEE Intelligent Transportation Systems Magazine}, 
  title={Making Bertha Drive—An Autonomous Journey on a Historic Route}, 
  year={2014},
  volume={6},
  number={2},
  pages={8-20}
  }

@ARTICLE{Broggi:2013,
  author={Broggi, Alberto and Buzzoni, Michele and Debattisti, Stefano and Grisleri, Paolo and Laghi, Maria Chiara and Medici, Paolo and Versari, Pietro},
  journal={IEEE Transactions on Intelligent Transportation Systems}, 
  title={Extensive Tests of Autonomous Driving Technologies}, 
  year={2013},
  volume={14},
  number={3},
  pages={1403-1415}
  }

@inproceedings{Madas:2013,
author = {Madås, David and Nosratinia, Mohsen and Keshavarz, Mansour and Sundströ, Peter and Philippsen, Roland and Eidehall, Andreas and Dahlén, Karl-Magnus},
year = {2013},
month = {06},
pages = {},
title = {On Path Planning Methods for Automotive Collision Avoidance},
}

@article{Isermann:2008,
author = { R.   Isermann  and  M.   Schorn  and  U.   Stählin },
title = {Anticollision system PRORETA with automatic braking and steering},
journal = {Vehicle System Dynamics},
volume = {46},
number = {sup1},
pages = {683-694},
year  = {2008},
publisher = {Taylor & Francis},


}

@inproceedings{Ploeg:2022,
  author={van der Ploeg, C and Smit, R. and Teerhuis, A. and Silvas, E.},
  booktitle={2022 IEEE 25th International Conference on Intelligent Transportation Systems (ITSC)}, 
  title={Long Horizon Risk-Averse Motion Planning: A Model-Predictive Approach}, 
  year={2022},
  volume={},
  number={},
  pages={1-8}}

@inproceedings{Smit:2022,
  author={Smit, R. and van der Ploeg, C. and Teerhuis, A. and Silvas, E.},
  booktitle={2022 IEEE 25th International Conference on Intelligent Transportation Systems (ITSC)}, 
  title={Informed sampling-based trajectory planner for automated driving in dynamic urban environments}, 
  year={2022},
  volume={},
  number={},
  pages={1-8}}

@INPROCEEDINGS{Stiller2010,  author={Ziegler, Julius and Stiller, Christoph},  booktitle={2010 IEEE Intelligent Vehicles Symposium},   title={Fast collision checking for intelligent vehicle motion planning},   year={2010},  volume={},  number={},  pages={518-522}}

@misc{Naoyuki:17,
	author={Tsushima {Naoyuki} and Abukawa {Masahiro}},
	title={Collision risk calculation device, collision risk display device, and vehicle body control device},
	assignee={ Mitsubishi Electric},
	type={US Patent},
	note={20170186319 A1},
	year={2017},
}

@misc{Sven:17,
	author={Heinig {Sven} and Hoever {Norbert} and Magyar {Peter} and Ottenhues {Thomas} and Wassmuth {Joachim}},
	title={Path planning},
	assignee={Hella KGaA Hueck \& Co.},
	type={US Patent},
	note={20080255728 A1},
	year={2008},
}
\end{document}